\documentclass[3p,times,review]{elsarticle}
\usepackage{color}
\usepackage{times}
\usepackage{booktabs}
\usepackage[pdfpagelabels=true,plainpages=false,colorlinks=true,linkcolor=blue,citecolor=blue,urlcolor=blue]{hyperref}
\usepackage[latin1]{inputenc}
\usepackage{algorithm}
\usepackage{algorithmic}
\usepackage{amsfonts}
\usepackage{amsmath}
\usepackage{amsthm}
\usepackage{tikz}
\usetikzlibrary{shapes.geometric, arrows}
\newcommand{\be}{\begin{equation}}
\newcommand{\ee}{\end{equation}}
\newcommand{\bea}{\begin{eqnarray}}
\numberwithin{equation}{section}
\newcommand{\eea}{\end{eqnarray}}
\newtheorem{thm}{Theorem}[section]

\newtheorem{remark}{Remark}
\makeatletter
\newcommand{\vast}{\bBigg@{4}}
\newcommand{\Vast}{\bBigg@{5}}
\makeatother
\usepackage{graphicx,epsfig}
\usepackage{enumerate}
\usepackage{setspace}
\usepackage{amssymb}
\usepackage{amsmath}
\usepackage{amsthm}
\usepackage{appendix}

\usepackage{subfig}
\theoremstyle{definition}

\newtheorem{dfn}{Definition}[section]
\newtheorem{note}{Note}[dfn]

\newtheorem{lem}{Lemma}[section]
\newtheorem{example}{Example}[section]

\usepackage{amssymb}
 \usepackage{amsthm}
\usepackage{ecrc}

\volume{00}

\firstpage{1}

\runauth{}

\usepackage{amssymb}

\makeatletter
\def\ps@pprintTitle{%
	\let\@oddhead\@empty
	\let\@evenhead\@empty
	\def\@oddfoot{}%
	\let\@evenfoot\@oddfoot}
\makeatother
\begin{document}
\begin{frontmatter}

\author[add1]{Shubham Kumar}
\address[add1]{Discipline of Mathematics, PDPM IIITDM Jabalpur, Jabalpur-482005, Madhya Pradesh, India}
 \ead{kumarshubham3394@gmail.com}


\author[add1]{Nihar Kumar Mahato\corref{cor1}}
\ead{nihar@iiitdmj.ac.in}
\cortext[cor1]{Corresponding author}

\author[add3]{Debdas Ghosh}
\address[add3]{Department of Mathematics, IIT(BHU), India}
\ead{debdas.mat@iitbhu.ac.in}
\title{Robust Optimization Approach for Solving Uncertain Multiobjective Optimization Problems Using the Projected Gradient Method}

\begin{abstract}
	Numerous real-world applications of uncertain multiobjective optimization problems (UMOPs) can be found in science, engineering, business, and management. To handle the solution of uncertain optimization problems, robust optimization is a relatively new field. An extended version of the projected gradient method (PGM) for deterministic smooth multiobjective optimization problem (MOP) is presented in the current study as a PGM for UMOP. An objective-wise worst-case cost (OWWC) type robust counterpart is considered, and the PGM is used to solve a UMOP by using OWWC. A projected gradient descent algorithm is created using theoretical findings. It is demonstrated that the projected gradient descent algorithm's generated sequence converges to the robust counterpart's weak Pareto optimal solution, which will be the robust weak Pareto optimal solution for UMOP. Under a few reasonable presumptions, the projected gradient descent algorithm's full convergent behavior is also justified.  Finally, numerical tests are presented to validate the proposed method.

\end{abstract}


\end{frontmatter}

\label{}
\section{Introduction}
\label{intro}
Optimization plays a crucial role in various fields like engineering, industry, medicine, and business, as well as across different scientific disciplines. Real-world optimization problems often involve dealing with uncertain data. This uncertainty can arise from inherent randomness or errors associated with the data. For instance, when estimating the demand for a product in inventory management, inaccuracies may occur due to a lack of knowledge about certain parameters of the mathematical model. Handling parameter uncertainty poses a significant challenge in solving optimization problems. In recent times, addressing uncertainty in optimization problems has garnered considerable attention within the mathematical programming community. In many real-world optimization problems, there is not just one clear objective function. Instead, the goals depend on various factors because multiple decision-makers are involved, each with their own optimization criteria. In real-world scenarios, optimization often faces the challenge of dealing with both uncertainty and multiple objectives. For instance, when deciding how much inventory to keep for a product, there is uncertainty about the demand, and different stakeholders might have different goals, like minimizing costs or maximizing customer satisfaction. Handling these uncertainties while considering various objectives simultaneously is a complex task. Researchers have studied two main areas: one deals with uncertain data using stochastic and robust optimization, while the other looks at different goals through multiobjective decision-making and optimization.
\par In the last decades, many researchers are focusing on developing methods to tackle UMOPs effectively \cite{ide2016robustness, koushki2022lr, kruger2023point, schmidt2019min, schobel2021price, shavazipour2021multi, zhou2019decision}. As one can observe that most of the  MOPs suffer from uncertain data. To understand the uncertainty in MOPs one real word example can be consider: Consider a manufacturing company that produces electronic devices. The company aims to optimize the design of a new smartphone model, considering multiple conflicting objectives such as maximizing battery life, minimizing production cost, and maximizing processing speed. However, due to uncertainties in factors such as market demand, component availability, and technological advancements, the company faces uncertainty in several parameters of the optimization problem. Continuing with the smartphone manufacturing company, consideration is given to a scenario where the uncertainties are bounded within a finite set.\\
\textbf{Objective Functions:}\\
$~~~~~\bullet$ \textbf{Maximize Battery Life (Objective 1):} The company aims to maximize the battery life of the smartphone.\\ 
	~~$~~~~~\bullet$~\textbf{Minimize Production Cost (Objective 2):} Minimize the overall production cost of the smartphone.\\
	~~$~~~~~\bullet$~ \textbf{Maximize Processing Speed (Objective 3):} Maximize the processing speed of the smartphone.\\
\textbf{Constraints:}\\
$~~~~~\bullet$ \textbf{Battery Capacity Constraint:} The battery capacity must be within a specified range.\\
~~	$~~~~~\bullet$ \textbf{Component Availability Constraint:} Availability of key components within a certain range.\\
~~$~~~~~\bullet$ \textbf{Technological Constraints:} The design must meet size limitations and compatibility requirements.\\
~~	$~~~~~\bullet$ \textbf{Market Demand Constraint: } The production volume must meet a minimum threshold to meet expected market $~~~~~~~~$demand.\\
\textbf{Uncertainty:}\\
$~~~~~\bullet$ \textbf{Demand Uncertainty:} Market demand varies within a known range due to factors like seasonal trends, economic $~~~~~~~~$conditions, and marketing efforts.\\
	$~~~~~\bullet$ \textbf{Component Cost Uncertainty:} The prices of key components fluctuate within known upper and lower bounds $~~~~~~~~$due to market conditions and supplier contracts.\\
	$~~~~~\bullet$ \textbf{Battery Life Uncertainty:} The actual battery life of the smartphone varies within a specified range due to\\ $~~~~~~~~$manufacturing variations and performance characteristics.\\
By using robust optimization approach, for given a finite uncertainty sets for market demand, component costs, and battery life, the company aims to find robust solutions that perform well across all possible scenarios within these uncertainty sets.
This example illustrates how uncertain constrained MOPs can be addressed under finite uncertainty sets, providing practical insights for decision-making in complex and uncertain environments. Apart from this example, such type of many problems can be found in the real world. Therefore, a generalized mathematical model for such type of problems can be consider as follows:
\begin{equation}\label{1.1*}
	P(U)=\{P(x,\xi_i):\xi_i\in U\},
\end{equation} where $U=\{\xi_i:i\in \bar\Lambda=\{1,2,\ldots,p\}\}$ is a set of uncertain parameters and $P(x,\xi_i)$ is given by
$$P(x,\xi_i):~~~~~~\min_{x\in D} h(x,\xi_i).$$ For $P(U),$  $h:\mathbb{R}^n\times U\to \mathbb{R}^m$ such that $h(x,\xi_i)=(h_1(x,\xi_i),h_2(x,\xi_i),\ldots,h_m(x,\xi_i)),$ where $h_j:\mathbb{R}^n\times U\to\mathbb{R},$ $j=1,2,\ldots,m.$
To handle the solution of $P(U),$ the concept of robust optimization and projected gradient method which is developed for deterministic MOPs are used in this article. 
\par In single objective optimization, the objective is to find a single solution that optimizes a specific function, such as maximizing profit or minimizing cost. This solution typically represents a single point in the solution space. On the other hand, deterministic MOP contains more than one objective functions. Because different objectives often conflict with each other, it is impossible to find a solution that perfectly optimizes all of them at once.  Instead, a set of solutions, known as efficient or Pareto optimal solutions, are available. An efficient or Pareto optimal solution is one where improving any objective function results in a decline in at least one other objective function. The value of the objective functions at an efficient solution is termed a nondominated point. The collection of all nondominated points and the set of all efficient solutions are referred to as the nondominated and efficient sets, respectively (for further details, see \cite{miettinen1999nonlinear, ehrgott2005multicriteria}). 
Various methods have been introduced to find the Pareto optimal point of a MOP: Weighted sum methood \cite{eckenrode1965weighting}; Adaptive weighted sum method \cite{koski1985defectiveness}; An $\epsilon-$constraint method \cite{haimes1971bicriterion}; Normal boundary intersection method \cite{das1998normal}; Normal constraint method \cite{messac2004normal}; Direct search domain method \cite{erfani2011directed};  Ideal cone method  \cite{ghosh2014new}. One can find more about the solution of MOP in \cite{marler2004survey, gunantara2018review, miettinen2008introduction} and references therein.
In addition to aforementioned techniques, classical derivative methods have been extended to solve MOPs: Steepest descent method \cite{fliege2000steepest,drummond2005steepest}, initially proposed for multicriteria optimization and then for vector optimization, this method converges to a critical point of the objective function; Projected gradient method \cite{drummond2004projected} developed for constrained MOPs, this method converges to critical and weak efficient points for convex and nonconvex MOPs, respectively. Other classical derivative methods \big(Projected gradient method \cite{cruz2011convergence,fazzio2019convergence,fukuda2013inexact,fukuda2011convergence,zhao2021projected, zhao2018convergence}, Newton's method \cite{fliege2009newton}, quasi-Newton method \cite{lai2020q,ansary2015modified,povalej2014quasi,mahdavi2020superlinearly,morovati2018quasi,qu2011quasi}, conjugate gradient method \cite{gonccalves2020extension,lucambio2018nonlinear}, proximal gradients \cite{bonnel2005proximal,ceng2010hybrid}\big) are also extended from scalar optimization to MOPs.
\par In uncertain optimization, different methods have been explored, including stochastic optimization and robust optimization. Stochastic optimization, as explained by Birge and Louveaux \cite{birge2011introduction}, relies on probabilistic information about uncertainties. It involves concepts like expected value, chance constraints, and risk measures such as conditional value-at-risk, which are also used in multiobjective problems. On the other hand, robust optimization does not rely on probabilistic information. Instead, it assumes that uncertain parameters come from a defined uncertainty set. These uncertain parameters are often referred to as scenarios. 

To handle the uncertainty in uncertain single objective optimization problems (USOPs) by robust optimization approach, there are different ideas about what makes a solution desirable. One well-known concept is called minmax (or strict) robustness, first talked about by Soyster \cite{soyster1973convex} and extensively studied by Ben-Tal and Nemirovski \cite{ben1998robust,ben1999robust}. In minmax robustness, a solution is considered robust if it is feasible for every possible scenario and if it minimizes the objective function in the worst-case scenario. This means it performs well even under the worst conditions, but it might be overly cautious because it has to work for every possible scenario. Instead of just focusing on minimizing the worst-case scenario, other objective functions like absolute or relative regret can also be used. These alternatives aim to capture how much a solution deviates from the best possible outcome. Because a robust solution needs to work for every scenario, these concepts are sometimes criticized for being too cautious. Recently, there have been some new ideas introduced to address this issue \big(see in Fischetti and Monaci \cite{fischetti2009light}; Schöbel \cite{schobel2014generalized}; Ben-Tal et al. \cite{ben2004adjustable}; Ben-Tal et al. \cite{ben2010soft}; Liebchen et al. \cite{liebchen2009concept}; Erera et al. \cite{erera2009robust}; Goerigk and Schöbel \cite{goerigk2014recovery}; Goerigk and Schöbel \cite{goerigk2016algorithm}\big).
\par Dealing with uncertainties in MOPs, RO is a relatively new area of research. Early research in this field did not focus on the typical concepts of robustness. Instead, most of the work centered around an idea introduced by Branke \cite{branke1998creating}, which was initially developed for USOPs. In this concept, the objective function is replaced by its average value within a specific neighborhood around the chosen point. Based on this concept, Deb and Gupta \cite{deb2013multi} introduced, two approaches to handling uncertainties in MOPs. In first approach, the objective vector is changed to a vector that contains the average values of each original component. Therefore, robust solutions to the original problem are efficient solutions to the modified problem. In second approach, mean functions are included in the constraints, which makes sure the objective components do not deviate from their average values by a distance greater than a set limit. Because users can control the desired level of robustness through this threshold, the latter approach is thought to be more practical. 
The concept of robustness for USOPs is expanded for UMOPs by an approach that Doolittle et al. introduced in \cite{doolittle2018robust}. Kuroiwa and Lee also introduced a similar strategy in \cite{kuroiwa2012robust}. Doolittle et al. addressed UMOPs by incorporating additional variables, following the method used by Ben-Tal and Nemirovski \cite{ben1998robust} for USOPs. A worst-case scenario for one objective function is represented by each variable. Solutions that perform well in this revised problem are thought to be robust. These new variables produce a new objective vector. Avigad and Branke \cite{avigad2008embedded} proposed a unique extension of minimax robustness for UMOPs. Instead of applying the worst-case scenario to individual components, this approach applied it to the entire objective vector. This created a deterministic MOP over the uncertainty set, where the worst case represented a set of scenarios. Avigad and Branke \cite{avigad2008embedded} developed an evolutionary algorithm to find robust solutions, meaning solutions where the worst-case set could not be dominated by another worst-case set. In simpler terms, the Avigad and Branke's approch aimed to find solutions not contained in a set resulting from a different worst-case scenario, with specific conditions fixed. Ehrgott et al.~\cite{ehrgott2014minmax} extended the minimax robustness concept in a general way for multiobjective optimization from the single objective optimization. 
Ehrgott et al.~\cite{ehrgott2014minmax} solved the UMOP with the help of the scalarization method (e.g., weighted sum approach and $\epsilon-$constraints approach) by using a general extension of the minimax concept of the robustness from USOP to UMOP. Using RO approach, Ehrgott et al. \cite{ehrgott2014minmax} changed the UMOP into a deterministic MOP with the help of a minimax type robust counterpart and objective-wise worst-case cost type robust counterpart. Ehrgott et al.~\cite{ehrgott2014minmax} solved robust counterparts of UMOPs with the help of scalarization techniques like the weighted sum approach and the $\epsilon-$constraint approach. There are some drawbacks to using scalarization techniques to solve the robust counterparts of UMOPs, such as the pre-specification of weights, restrictions, or function importance that is unknown beforehand. 
To overcome these difficulties, in recent years, there has been a lot of work on developing robustness concepts for dealing with uncertainty in multiobjective optimization problems \big(see in  Ide and Schöbel \cite{ide2016robustness}; Kaushki et al. \cite{koushki2022lr}; Krugelet al. \cite{kruger2023point}; Schmidt et al. \cite{schmidt2019min}; Schöbel and Kangas \cite{schobel2021price}; Shavazipour et al. \cite{shavazipour2021multi}; Shavazipour et al. \cite{shavazipour2021multia};
Kangas and Miettinen \cite{zhou2019decision}; Bokrantz and Fredriksson \cite{bokrantz2017necessary}\big). Ide and Schöbel \cite{ide2016robustness} conducted a survey and analysis of different concepts related to robustness in uncertain multiobjective optimization. Wiecek and Dranichak \cite{wiecek2016robust} investigated methods for robust multiobjective optimization to help make decisions when facing uncertainty and conflict.
Also, the concept of minmax robust efficient has attracted a lot of attention in various fields. Some studies take a set-valued perspective, others focus on specific applications like portfolio selection or linear uncertainty cases \cite{fliege2014robust, wiecek2016robust}. There are also different approaches for dealing with uncertainty in constraints or in problems like shortest path  Kalantari et al. \cite{kalantari2016multi}. Bokrantz and Fredriksson \cite{bokrantz2017necessary}, along with Ehrgott et al. \cite{ehrgott2014minmax}, as well as Schmidt et al. \cite{ schmidt2019min}, propose scalarization methods. 
Since minmax robust efficient solutions tend to be overly cautious, other concepts have been developed. These include highly robust efficiency, lightly robust efficiency, and regret robustness. Scenario-based versions of these concepts have also been explored. In recent years, analyzing solutions to uncertain multiobjective problems has gained significant attention. Schöbel and Zhou-Kangas \cite{schobel2021price} highlight the challenge of defining robust efficient solutions and propose an approach to compare these solutions, and introduce insights into different solution sets, such as nominal efficient, minmax robust efficient, and lightly robust efficient solutions. Also provides a measure to evaluate the robustness of individual solutions, offering strategies for decision-making in uncertain multiobjective optimization. 
Koushki et al. \cite{koushki2022lr} present LR-NIMBUS, an interactive algorithm designed to aid decision-makers in finding lightly robust efficient solutions for uncertain multiobjective optimization problems.
Krügera et al. \cite{kruger2023point} introduce a novel concept termed ``the point based robustness gap''. This concept addresses the challenge of uncertain multiobjective optimization problems, and propose a measure of the distance between the robust Pareto set and the Pareto sets of scenarios, termed as the robustness gap. 
\par Apart form the above techniques some iterative methods is also developed for uncertain unconstrained multiobjective problems. In this direction of research, Kumar et al. \big(\cite{shubham2023newton, kumar2023quasi, kumar2024quasi}\big) considered an uncertain unconstrained multiobjective optimization problem under finite uncertainty set. To tackle the solution of this problem Kumar et al. used objective-wise worst-case-type robust counterpart. The objective-wise worst-case-type robust counterpart is solved by Newton's method \cite{shubham2023newton}, quasi-Newton method \cite{kumar2023quasi}, modified quasi-Newton method \cite{kumar2024quasi} to find the solution of UMOP. Since these methods are applicable for uncertain unconstrained multiobjective optimization problems and not applicable for uncertain constrained multiobjective optimization problems. 
In this article, a modified version of projected gradient method is developed to find the solution of uncertain constrained multiobjective optimization problem given in equation (\ref{1.1*}). Drummond and Iusem \cite{drummond2004projected} have been extended the classical projected gradient method to solve vector optimization problems. Drummond and Iusem \cite{drummond2004projected} have been presented that the projected gradient method converges to critical and weak efficient points for nonconvex and convex MOPs, respectively. Further development related to projected gradient methods can be seen in the literature \cite{cruz2011convergence,fazzio2019convergence,fukuda2013inexact,fukuda2011convergence,zhao2021projected, zhao2018convergence} and references therin. 
To the best knowledge of the authors, there is no projected gradient method developed for uncertain constrained multiobjective optimization problems.
\par The article's organization is as follows: Important results, fundamental definitions, and theorems essential to the problem are presented in Section \ref{s2}. Section \ref{s3} begins by providing the Pareto optimality condition for OWWC. Subsequently, in Subsection \ref{ss3.1}, a search direction subproblem is solved to find the projected gradient descent direction. Following that, in Subsection \ref{ss3.2}, an Armijo-type inexact line search technique is established to determine an appropriate step length size, ensuring that the function value decreases along the projected gradient descent direction. Subsection \ref{algo1} presents the projected gradient descent algorithm for OWWC. By using this algorithm, a sequence is generated, and its convergence to a critical point is proven in Subsection \ref{ss3.3}. In Subsection \ref{s341}, under certain assumptions, the full convergence of the sequence generated by the projected gradient descent algorithm converges to a Pareto optimum is shown. In Section \ref{secnum}, some numerical tests are presented to validate the projected gradient method. Section \ref{s4} concludes the article with some remarks on the projected gradient descent method.
\section{General concepts in uncertain multiobjective optimization and deterministic multiobjective optimization}\label{s2}
Beginning with some notations. The set of real numbers denoted as $\mathbb{R}$, the set of non-negative real numbers denotes as $\mathbb{R}_{\geq} = \{x \in \mathbb{R} : x \geq 0\}$, the set of positive real numbers denotes as $\mathbb{R}_{>} = \{x \in \mathbb{R} : x > 0\}.$ 
$\mathbb{R}^{n} = \mathbb{R} \times \cdots \times \mathbb{R} ~(n~ \text{times}),$ $\mathbb{R}^{n}_{\geqq} = \mathbb{R}_{\geq} \times \cdots \times \mathbb{R}_{\geq}$ ($n$ times), $\mathbb{R}^{n}_{\geq}=\mathbb{R}^{n}_{\geqq}\setminus\{0\},$ $\mathbb{R}^{n}_{>} = \mathbb{R}_{>} \times \cdots \times \mathbb{R}_{>}$ ($n$ times).
In the context of $s, q \in \mathbb{R}^n$, the notation $s \leq q$ denotes that the component of $s$ is related to the corresponding component of $q$ such that $s_i \leq q_i,$ for all $i = 1, 2,\ldots, n$.
Similarly, for any $s, q \in
\mathbb{R}^n$: $ s\geq q\Longleftrightarrow	s - q \in \mathbb{R}^n_{\geq}$, which is equivalent to $s_i - q_i \geq 0$ for each $i$,
$s > q \Longleftrightarrow	 s - q \in \mathbb{R}^n_{>}$, which is equivalent to $s_i - q_i > 0$ for each $i$.
Lastly, the indexed sets are denoted as $\bar\Lambda = \{1, 2,\ldots, p\}$ and $\Lambda = \{1, 2,\ldots, m\}$, containing $p$ and $m$ elements, respectively.
\par In single objective optimization, the objective is to find a single solution that optimizes a specific function, such as maximizing profit or minimizing cost. This solution typically represents a single point in the solution space.
On the other hand, multiobjective optimization involves optimizing multiple conflicting objectives simultaneously. These objectives may contradict each other, making it challenging to optimize one without compromising others. In multiobjective optimization, the solution space comprises a set of solutions known as the Pareto front. Here, no solution is superior to another in all objectives. Instead, each solution represents a trade-off between different objectives. The aim is to identify a set of solutions that collectively strike a balance across all objectives.
\par A MOP can be considered as 
$P:~\displaystyle\min_{x\in D \subset \mathbb{R}^n}F(x),$ where $F:D\to\mathbb{R}^m$ and $F(x)=(F_1(x),F_2(x),$\ldots$,F_m(x)).$ If $m=1,$ then the problem $P$ is called single objective optimization problem. A feasible point $x^*\in D$ is said to be optimal solution for $F$ (for $m=1$) if $F(x^*)\leq F(x)$ for all $x\in D.$ On the other hand, the solution related to MOP (i.e., $m>1$) defined as follows.
\begin{dfn} Given a point $x^* \in D$,
\begin{enumerate}
	\item[(i)] $x^*$ is said to be an efficient solution or Pareto optimal solution of $F$ if and only if there does not exist $x \in D$ such that $F(x) \leq F(x^*)$ and $F(x) \neq F(x^*)$;
	\item[(ii)] $x^*$ is said to be weakly efficient or weak Pareto optimal solution of $F$ if and only if there does not exist $x \in D$ such that $F(x) < F(x^*)$;
	\item[(iii)] $x^*$ is said to be locally efficient or locally weakly efficient for $F$ if and only if there exists a neighborhood $X \subset D$ such that the point $x^*$ is efficient or weakly efficient for $F$ restricted to $X$.
\end{enumerate}
\end{dfn}
If $x^*$ is an efficient (weakly efficient) solution, then $F(x^*)$ is called a non-dominated (weakly non-dominated) point, and the set of efficient solutions and the set of non-dominated points (Pareto front) are called the efficient set and the non-dominated set, respectively. For further details of these definitions, see references \cite{cruz2011convergence,fazzio2019convergence,fukuda2013inexact,fukuda2011convergence,zhao2021projected, zhao2018convergence}.\\  
\par In uncertain multiobjective optimization, input data that are uncertain affect how the optimization problem is formulated. This uncertainty is represented as a set $U\subset \mathbb{R}^k,$ which includes all the scenarios (or realization) of the uncertain data. Recall the uncertain constrained multiobjective optimization problem defined in the introduction section in equation (\ref{1.1*}) for further consideration.
\begin{equation*}\label{1.1}
	P(U)=\{P(x,\xi_i):\xi_i\in U\},
\end{equation*} where $U=\{\xi_i:i\in \bar\Lambda\}$ is a set of uncertain parameters and $P(x,\xi_i)$ is given by
$$P(x,\xi_i):~~~~~~\min_{x\in D} h(x,\xi_i).$$ For $P(U),$  $h:\mathbb{R}^n\times U\to \mathbb{R}^m$ such that $h(x,\xi_i)=(h_1(x,\xi_i),h_2(x,\xi_i),\ldots,h_m(x,\xi_i)),$ where $h_j:\mathbb{R}^n\times U\to\mathbb{R},$ $j\in \Lambda.$ Additionally, for any given $x\in D$ to $P(U)$, the set of images of $x$ under all scenario is given by $h_U(x)=\{h(x,\xi_i):\xi_i\in U\}.$
\par In this situation, it is crucial to note that the uncertainty specifically pertains to the objective function and not the constraints. This distinction arises because a solution to the robust problem is deemed feasible only under the standard minimax robustness definition if it remains feasible across all scenarios. This particular type of uncertainty is referred to as parameter uncertainty. In the constrained uncertain MOP $P(U)$ under consideration, uncertainty exists solely in the objective functions, not in the decision variables. Therefore, we are exclusively dealing with parameter uncertainty here, with no variability in the decision variables themselves. However, if $D$ is also considered uncertain, every robust solution $x$ to the uncertain problem must satisfy the constraints $x \in \bigcap_{i \in \bar{\Lambda}} D(\xi_i)$, where $D(\xi_i)$ denotes the feasible set corresponding to $\xi_i$. This kind of uncertainty is known as decision uncertainty. For further details, see Schöbel \cite{schobel2014generalized} and Ehrgott et al. \cite{ehrgott2014minmax}.
\par It can be observed that $P(U)$ is a set-valued optimization problem. So the concept of optimality of $P(U)$ can be adopted from the set-valued optimization. There is a need to compare the sets that represent the various outcomes to find the best solution to a set-based problem. The so-called set approach, as described in Eichfelder and Jahn \cite{eichfelder2011vector}, Ha and Jahn \cite{jahn2011new}, and Kuroiwa \cite{kuroiwa1998natural}, compares sets using order relations. Here, an order relation $\preceq$ is used to compare sets by means of a given closed convex pointed solid cone $C\subset\mathbb{R}^m.$ The relation $\preceq_C,$ $\preceq_{C\symbol{92}\{0\}},$ and $\preceq_{int\{C\}}$ are defined as follows for every sets $V,~ Q\subset \mathbb{R}^m,$
\begin{eqnarray*}
V\preceq_C Q \iff V\subset \{Q\}-C,\\
V\preceq_{C\symbol{92}\{0\}} Q \iff V\subset \{Q\}-C\symbol{92}\{0\},\\\mbox{and}~~~~
V\preceq_{int\{C\}} Q \iff V\subset \{Q\}-int\{C\}.
\end{eqnarray*}
Equivalently, these set relations can be represented as 
\begin{eqnarray}
V[\preceq_{int\{C\}}\symbol{92}{\preceq_{C\symbol{92}\{0\}}}\symbol{92}\preceq_C ]\iff \text{for all}~ v\in V~ \exists q\in Q ~\text{such that}~[v<\symbol{92}\leq\symbol{92}\leqq q].
\end{eqnarray} 
Throughout the article, a closed convex pointed solid cone $C=\mathbb{R}^m_\geq$ will be considered. On behalf of this cone the concept of robust efficiency can be defined in the following way. \\
\begin{dfn}(Ehrgott et al. \cite{ehrgott2014minmax}) Given a UMOP $P(U),$ a feasible point $x^*\in D$ is said to be 
\begin{itemize}
	\item robust efficient (RE) if there is no  $x\in D-\{x^*\}$ such that $h(x;U)\subset h(x^*;U)-\mathbb{R}^k_\ge$ ,
	\item robust weakly efficient (RWE) if there is no  $x\in D-\{x^*\}$ such that $h(x;U)\subset h(x^*;U)-\mathbb{R}^k_>$,
	\item robust strictly efficient (RSE) if there is no  $x\in D-\{x^*\}$ such that $h(x;U)\subset h(x^*;U)-\mathbb{R}^k_\geqq$.
\end{itemize}
\end{dfn}
In \cite{ehrgott2014minmax}, the authors introduce an extension of the concept of minimax robustness for USOPs, as introduced by Soyster \cite{soyster1973convex} and studied, e.g., by Ben-Tal and Nemirovski \cite{ben1998robust}. Further, a given UMOP $P(U),$ the concept of minimax robustness is presented by Ehrgott et al. \cite{ehrgott2014minmax}. The motive of this robustness concept is to search for solutions that minimize the worst case that occur, i.e., solution to the problem
\begin{eqnarray}\label{minimax}
\displaystyle 	\min_{x\in D}\sup_{\xi_i\in U}h(x,\xi_i),~\text{where}~ h:D\times U\to \mathbb{R}^m.
\end{eqnarray}	
It is pointed out that the definition of a worst-case is not clear anymore due to the lack of total order on $\mathbb{R}^m.$ Thus an extension of the concept of minimax robustness to UMOP was not directly possible. Therefore, Ehrgott et al. \cite{ehrgott2014minmax} presented an extension of minimax robustness to MOPs, namely the concept of robust efficiency. This concept is presented above in the form of RE solution, RWE solution, and RSE solution. 
\par Ehrgott et al. \cite{ehrgott2014minmax} developed algorithms for computing minimax RE solution. Firstly, a well-known weighted sum scalarization method for deterministic MOP is extended to a weighted sum scalarization method for finding minimax RE solutions to UMOPs.
One more approach is presented by Ehrgott et al. \cite{ehrgott2014minmax} for finding the minimax robust efficient solution to $P(U).$ In this approach the objective-wise worst-case (OWWC) method is considered. This approach consists of the same deterministic multiobjective optimization problem as Kuroiwa and Lee \cite{kuroiwa2012robust} considered for their concepts of multiobjective robustness, namely    
\begin{eqnarray}\label{owc}
OWC_{P(U)}:~~~~~~~\min_{x\in D} h^{owc}_U(x),~\text{where}~ h^{owc}_U(x)=(\sup_{\xi_i\in U}h_1(x,\xi_i),\sup_{\xi_i\in U}h_2(x,\xi_i),\ldots,\sup_{\xi_i\in U}h_m(x,\xi_i))^T.
\end{eqnarray}
It is also presented by Ehrgott et al. \cite{ehrgott2014minmax} that the solution to $OWC_{P(U)}$ will be the solution for $P(U).$
\begin{thm}\label{t1}( Ehrgott et al. \cite{ehrgott2014minmax}) Let $P(U)$ be an UMOP. 
\begin{enumerate}[(a)]
	\item If $x^*\in X$ is a strictly efficient solution to $OWC_{P(U)}$, then $x^*$ is the RSE solution for $P(U)$.
	\item If $\displaystyle\max_{\xi_i\in U} h_j(x,\xi_i)$  exist for all $j\in\Lambda$ and all $x \in X$ and $x^*$
	is a weakly efficient solution to $OWC_{P(U)}$, then $x^*$ is RWE solution for $P(U)$.
\end{enumerate}
\end{thm}
In addition to this result, Ehrgott et al. \cite{ehrgott2014minmax} developed the weighted sum scalarization method and $\epsilon$-constraints scalarization method for $OWC_{P(U)}$ to find the solution for $P(U)$. It has also been proven that the solutions obtained by both scalarization methods for $OWC_{P(U)}$ will be solutions for $P(U).$
\begin{remark} Obviously, computing $h^{owc}U(x)$ for any given $x$ is much easier than solving a MOP $\max\{h(x,\xi_i):\xi_i\in U\}$. Additionally, $OWC_{P(U)}$ contains only $m$ deterministic SOPs; therefore, it is considered a deterministic MOP. Hence, to find the solution of $P(U)$, $OWC_{P(U)}$ is considered instead of the minimax-type robust counterpart of $P(U)$.
\end{remark}
Theorem \ref{t1} shows that the solution of $OWC_{P(U)}$ will be the solution for $P(U).$ Additionally, the problem $OWC_{P(U)}$ is solved by scalarization methods (e.g., weighted sum method and $\epsilon-$constraint method). By pre-selecting some parameters and reformulating them as a deterministic SOP, scalarization methods, which are based on the scalarization technique, compute the efficient or weakly efficient solution. 
To solve deterministic multiobjective optimization problems, scalarization methods (such as the weighted sum method and $\epsilon-$constraints method) and other methods mentioned in the introduction encounter certain difficulties, as discussed in  \cite{eckenrode1965weighting,koski1985defectiveness,haimes1971bicriterion,messac2004normal,das1998normal,messac2004normal,ghosh2014new,erfani2011directed,marler2004survey, gunantara2018review, miettinen2008introduction}. For instance, when applying scalarization methods to uncertain multiobjective optimization problems, challenges arise. In the weighted sum method, selecting weights that keep the solution feasible across all scenarios is crucial. Often, chosen weights may unintentionally bias solutions towards unbounded regions, which is a drawback of this method \cite{ehrgott2014minmax,kuroiwa2012robust}. Similarly, applying the epsilon constraints method introduces difficulties in choosing epsilon values under uncertainty. Another drawback of these methods is that the parameters are not predetermined, so the responsibility lies with the modeler and the decision-maker to make those decisions.

Although, several new methods have been developed recently, as referenced in \cite{ide2016robustness,bokrantz2017necessary, koushki2022lr, kruger2023point,schmidt2019min,shavazipour2021multi,shavazipour2021multia,zhou2019decision}, to address these issues, scalarization methods still encounter challenges in nonconvex cases these methods are unable to find a complete Pareto front. In the case of unconstrained UMOP (i.e., $P(U)$ with $D=\mathbb{R}^n$),  descent methods like Newton's method \cite{shubham2023newton}, quasi-Newton method \cite{kumar2023quasi} and modified quasi-Newton method \cite{kumar2024quasi} have been developed for the $P(U)$ by using OWWC type robust counterpart of it. These descent methods overcome the difficulties of aforementioned techniques. However, these methods are not directly applicable to constrained UMOP $P(U).$ Therefore, we propose extending the projected pradient method, traditionally used for deterministic constrained multiobjective optimization problems, to handle constrained UMOP $P(U).$ 
\par It is also noted that existing methods are designed for smooth problems, while the development of the projected gradient method for non-smooth problems in the context of $OWC_{P(U)}$ is essential. To the best of our knowledge, no projected gradient method has been developed specifically to address the solution of constrained UMOP. Note that, to overcome the situation caused by the scalarization techniques and the other techniques given in \cite{ide2016robustness,bokrantz2017necessary, koushki2022lr, kruger2023point,schmidt2019min,shavazipour2021multi,shavazipour2021multia,zhou2019decision}, the projected gradient method is a good alternative, and no parameter information is needed for this method. This approach does not make use of predetermined weighting factors or any other kind of ranking or ordering data for the various objective functions.
\par A projected gradient method for $OWC_{P(U)}$ 
will be developed to find the solution of $P(U).$ It is assumed that each $h_j(x,\xi_i)$ is continuously differentiable with respect to $x$ for each $i$ and $j.$ Then, the problem $OWC_{P(U)}$  can be equivalently written as
\begin{eqnarray*}\label{main owc}
OWC_{P(U)}:~~~~~~~\min_{x\in D} h^{owc}_{U}(x),~\text{where}~ h^{owc}_{U}(x)=(\max_{\xi_i\in U}h_1(x,\xi_i),\max_{\xi_i\in U}h_2(x,\xi_i),\ldots,\max_{\xi_i\in U}h_m(x,\xi_i))^T,
\end{eqnarray*}
or
\begin{eqnarray}\label{main owc}
OWC_{P(U)}:~~~~~~~\min_{x\in D} H(x),~\text{where}~  H(x)=(H_1(x),H_2(x),\ldots,H_m(x))^T~\text{and}~H_j(x)=\displaystyle\max_{\xi_i\in U}h_j(x,\xi_i),~j\in\Lambda.
\end{eqnarray}
Since $U$ is a finite set and $h_j(x,\xi_i)$ is continuously differentiable for each $x$ and $\xi_i,$ then $\displaystyle\sup_{\xi_i\in U}h_j(x,\xi_i)=\displaystyle\max_{\xi_i\in U}h_j(x,\xi_i)$ due to the compactness of $U.$
\par One can observe that $OWC_{P(U)}$ is a non-smooth deterministic MOP.
The necessary conditions of Pareto optimality for $OWC_{P(U)}$ are now defined.
\begin{dfn}\label{criticalpoint}\cite{qu2014nonsmooth} A point $x^*\in D$ is said to be Pareto critical point for $OWC_{P(U)}$ if $I(Conv\{\displaystyle \cup_{j\in \Lambda} \partial H_j(x^*)\})\cap(-\mathbb{R}^{n}_{>})= \emptyset.$
\end{dfn}
By this definition, it can be noticed that $x^*$ is a Pareto critical point for $OWC_{P(U)}$ if and only if there is no $v\in D-\{x^*\}$ such that $\nabla h_j(x^*,\xi_i)^Tv<0$, for all $i\in I_j(x),$ $j\in \Lambda$, where $I_j(x)=\{i\in\bar{\Lambda}:h_{j}(x^*,\xi_i)=H_j(x^*)\}.$

\begin{note} By the definition of critical point for $OWC_{P(U)}$ and Theorem \ref{t1}, it can be concluded that the Pareto critical point for $OWC_{P(U)}$ will be the robust Pareto critical point for $P(U)$. Throughout the article, the Pareto critical point for $OWC_{P(U)}$ will be used as the robust Pareto critical point for $P(U)$.
\end{note}
\begin{dfn}
In $OWC_{P(U)},$ a vector $v$ is said to be descent direction for $H$  at $x$ if
$\nabla h_j(x,\xi_i)^Tv< 0, ~\forall j \in \Lambda$ and $i\in I_j(x).$
Also, $v$ is descent direction for $H(x)$ $\Longleftrightarrow$ there exists $\epsilon >0$ such that $H_j(x+\alpha v)< H_j(x) ~\forall~ j \in \Lambda~ \text{and} ~\alpha \in(0,\epsilon].$
\end{dfn}
\begin{thm}\label{directional derivative}\cite{dhara2011optimality} Let $H_j: \mathbb{R}^n \to \mathbb{R}^n$ denote a function such that $H_j(x)=\displaystyle\max_{i\in\bar{\Lambda}}h_j(x,\xi_i).$ Then:
\begin{enumerate}[(i)]
	\item In the direction $v$, the directional derivative of $H_j$ at $x$ is given by $H'_j(x,v)=\displaystyle\max_{i\in I_j(x)}\nabla h_j(x,\xi_i)^Tv,$ where $I_j(x) = \{i \in \bar{\Lambda} : h_j(x, \xi_i) = H_j(x)\}$.
	\item\label{ddii} The subdifferential of $H_j$ is  \[\partial H_j(x)= Conv \bigg(\bigcup_{i\in I_j(x)}\partial h_j(x,\xi_i)\bigg).\] Also, $x^*=\displaystyle arg\min_{x\in \mathbb{R}^n}H_j(x)$ $\Longleftrightarrow$ $0\in\partial H_j(x^*)).$ 
\end{enumerate}
\end{thm}
\section{Projected gradient method (PGM) for $OWC_{P(U)}$}\label{s3}
\par 
In this section, some theoretical results are developed, which are required to develop the PGM for the problem $OWC_{P(U)}$. To demonstrate the necessary and sufficient condition for the Pareto optimality of $OWC_{P(U)}$, Lemma \ref{lm1} is first provided.

\begin{lem}\label{lm1}
The point $x^*$ is a Pareto critical point for $H$ $\Longleftrightarrow$ $0\in Conv\bigg\{\displaystyle \bigcup_{j\in \Lambda}\partial H_j(x^*)\bigg\}$.
\end{lem}
\begin{proof}
Assume $x^*$ is Pareto critical point for $H$ then there must exists $d\in\bigcup_{j\in \Lambda}\partial H_j(x^*)$ such that 
\begin{equation}\label{ccc}
	v^Td\geq0,~\forall v\in\mathbb{R}^n.
\end{equation} On contrary, assume that $0\not\in Conv\bigg\{\displaystyle \bigcup_{j\in \Lambda}\partial H_j(x^*)\bigg\}.$ Since, $Conv\bigg\{\displaystyle \bigcup_{j\in \Lambda}\partial H_j(x^*)\bigg\}$ and $\{0\}$ are closed and convex sets then with the help of theorem of separation, there exists $v\in\mathbb{R}^n$ and $b\in \mathbb{R}$ such that 
$v^T0\geq b$ and $v^Td<b~\forall d\in Conv\bigg\{\displaystyle \bigcup_{j\in \Lambda}\partial H_j(x^*)\bigg\}.$ Jointly both inequality contradicts (\ref{ccc}). Hence,  $0\in Conv\bigg\{\displaystyle \bigcup_{j\in \Lambda}\partial H_j(x^*)\bigg\}$.\\ Conversely, it needs to be proven that if $0\in Conv\bigg\{\displaystyle \bigcup_{j\in \Lambda}\partial H_j(x^*)\bigg\}$, then $x^*$ is a Pareto critical point for $H$. For this purpose, define $\breve{H}(x)=\displaystyle\max_{j\in\Lambda} H_j(x)- H_j(x^*)$. Then, by item (\ref{ddii}) of Theorem \ref{directional derivative}, $\partial \breve{H}(x)= Conv \bigg\{\displaystyle\bigcup_{i\in\Lambda}\partial H_j(x)\bigg\}$. Hence, the assumption leads to $0\in Conv \bigg\{\displaystyle\bigcup_{i\in\Lambda}\partial H_j(x)\bigg\}$, implying $x^*=\arg\displaystyle \min_{x\in D}\breve{H}(x)$. On the contrary, if $x^*$ is not a Pareto critical point, then according to Definition \ref{criticalpoint}, there exists $s\in D$ such that $\nabla h_j(x^*,\xi_i)^Ts<0$, for all $i\in I_j(x^*)$, $j\in \Lambda$, i.e., $H'_j(x^*,s)<0$ for all $j$. Then there exists some $\eta>0$ sufficiently small such that $H_j(x^*+\eta s)<H_j(x^*)$ for all $j$ which implies $\breve{H}(x^*+\eta s)<0=\breve{H}(x^*)$ holds for some $(x^*+\eta s)\in D$. This contradicts the fact that $x^*=\arg\displaystyle \min_{x\in D}\breve{H}(x).$ As a consequence, the assumption that $x^*$ is not a Pareto critical point is incorrect, and $x^*$ is indeed a Pareto critical point for $H$.
\end{proof}
\begin{thm}\label{nce}
If  $h_j(x,\xi_i)$ is continuously differentiable and convex for each $j\in \Lambda$ and $\xi_i\in U$, then  $x^*\in D$ is a weak efficient solution solution for $OWC_{P(U)}$ if and only if
$$	0 \in conv \left( \displaystyle \cup_{j=1}^{m}\partial H_j(x^*) \right).$$
\end{thm}
\begin{proof}
Let  $x^*$ be a weak efficient solution solution for  $OWC_{P(U)}$. It must be shown that $0\in \text{Conv}{\displaystyle \cup_{j\in \Lambda}\partial H_j(x^*)}$. Since given function $h_j(x,\xi_i)$ is continuously differentiable and convex for each $j$ and $\xi_i \in U$, then $h_j(x,\xi_i)$ will be locally Lipschitz continuous for all $i\in \bar{\Lambda}$. Then $0\in Conv\{\displaystyle \cup_{j\in \Lambda}\partial H_j(x^*)\}$ (see Theorem 4.3 in \cite{makela2014nonsmooth} ).\\ Conversely, by assumption $0\in Conv\{\displaystyle \cup_{j\in \Lambda}\partial H_j(x^*)\}$ it is clear that $x^*$ is Pareto critical point. Then for atleast one $j^0$, it is established that $H'_{j^0}(x^*,d) \geq 0, ~\forall~ d\in D-\{x^*\}$. Now, by using the Definition \ref{directional derivative}, it follows that 
\begin{equation}\label{Equ1}
	\nabla h_{j^0}(x^*,\xi_i)^Td\geq 0, ~\forall~ d\in D,~ i\in I_{j^0}(x^*).
\end{equation}
By convexity of $H_j$ and $h_j(x,\xi_i)$, it is obtained that
\begin{center}
	$h_{j^0}(x,\xi_i)\geq h_{j^0}(x^*,\xi_i) + \nabla h_{j^0}(x^*,\xi_i)^T(x-x^*),~ \forall ~i\in I_{j^0}(x^*)$ and $x, ~x^* \in D.$
\end{center}
Since the last term of the latest inequality is positive by (\ref{Equ1}), it is established that
\begin{center}
	$ h_{j^0}(x,\xi_i)\geq h_{j^0}(x^*,\xi_i),~ \forall~ i\in I_{j^0}(x^*),$
\end{center}
and therefore
\begin{center}
	$ H_{j^0}(x)\geq H_{j^0}(x^*),~ \forall x\in D,$
\end{center}
i.e., $x^*$ is weak efficient solution.

\end{proof}

If $D=\mathbb{R}^n$ is considered in $OWC_{P(U)}$, then the iterative scheme related to the steepest descent algorithm can be given as $x^{k+1}=x^k+\alpha_k d^k$, $k \geq 0$, where $d^k=\underset{t\in \mathbb{R}^n}{\mathrm{argmin}} \left( \Theta_{x^k}(t)+\frac{1}{2}\|t\|^2 \right)$, $\Theta_{x^k}(t)= \displaystyle \max_{j\in \Lambda}\max_{i\in \bar{\Lambda}} \{ h_j(x^k,\xi_i)+\nabla h_j(x^k,\xi_i)^T t - H_j(x^k)\}$, and $\alpha_k$ is the step length size which satisfies the Armijo-type inexact line search. Such an algorithm is not directly applicable for $OWC_{P(U)}$ if $D \neq \mathbb{R}^n$, i.e., constrained optimization problems in which the decision variable is required to lie within a constraint set $D$. If the steepest descent algorithm is applied to solve $OWC_{P(U)}$, where $D \neq \mathbb{R}^n$, then the iterative point $x^k$ need not lie in $D$. In order to account for the presence of the constraints, the algorithms need to be changed. A simple change involves the introduction of a projection. The idea is as follows: if $x^k+\alpha_k d^k \in D$, then set $x^{k+1}=x^k+\alpha_k d^k$, as usual. On the other hand, if $x^k+\alpha_k d^k \notin D$, then project it back into $D$ before setting $x^{k+1}$. To develop the projected gradient method for $OWC_{P(U)}$, the subproblem to find the projected gradient direction is considered in the following subsection.

\subsection{\textbf{A subproblem to find projected gradient descent direction for $OWC_{P(U)}$}}\label{ss3.1}
The PGM for smooth MOP is now extended to $OWC_{P(U)}$. It is assumed that $D$ is a closed and convex subset of $\mathbb{R}^n$. For $x \in D$, the idea of the projected gradient descent direction, which is initially defined for smooth MOP by Drummond et al. \cite{drummond2004projected}, is extended. The extended version of search direction is given by
\begin{equation}\label{eq31}
t(x)=\underset{t\in D-\{x\}}{\mathrm{argmin}}  \left(\beta \Theta_{x}(t)+\tfrac{1}{2}\|t\|^2 \right),\\
\end{equation} where $\Theta_{x}(t)= \displaystyle \max_{j\in\Lambda}\max_{i\in \bar{\Lambda}}  \{ h_j(x,\xi_i)+\nabla h_j(x,\xi_i)^Tt-H_j(x)\}$ and $\beta>0$ is a parameter.
In the next step, the approach involves solving the subproblem \eqref{eq32}
\begin{equation}\label{eq32}
\text{S}_P(x):~~~\underset{t\in D-\{x\}}{\mathrm{min}}  \left(\beta \Theta_{x}(t)+\tfrac{1}{2}\|t\|^2 \right).\\
\end{equation}
Note that being the maximum of max-linear function, $\Theta_{x}(t)$ is convex, therefore the objective function $\Psi_x(t)=\beta\Theta_{x}(t)+\tfrac{1}{2}\|t\|^2 $ in the subproblem \eqref{eq32} is strongly convex. Then the subproblem \eqref{eq32} has a unique solution which will be the projected gradient descent direction for
$OWC_{P(U)}$. Hence projected gradient descent direction $t(x)$ for
$OWC_{P(U)}$ is always well-defined. If it is assumed that $t(x)$ is the solution of subproblem of \eqref{eq32}, then $0\in\partial\Psi_x(t(x)).$ Also, due the strong convexity of $\Psi_x(t),$ the first order optimality condition  for the subproblem \eqref{eq32} will be necessary and sufficient \big(for the optimality condition of a convex function, see Proposition 4.7.2 in \cite{bertsekas2003convex}\big). Therefore, there exists a feasible direction $t(x)$ and $u(x)\in\partial\Psi_x(t)$ such that 
\begin{equation}\label{eq33}
\langle u(x),t-t(x)\rangle\geq0,~\text{ for all }~t\in D-\{x\}.
\end{equation}
%
  
Now, the goal is to find the subgradient $u(x)$ of $\Psi_x(t(x))$  where $t(x)$ represents a feasible direction. By using the formula of subdifferential for the maximum of convex function (see in \cite{hiriart1993convex}), the subgradient of $\Psi_x(t)$ can be find in the following way:
%
	From the expression of $\Psi_x(t)$, there exist positive values $\lambda_{ij}(x) > 0$ for $j \in J_x \subseteq \Lambda$ and $i\in I_x \subseteq \bar\Lambda$ such that $\displaystyle\sum_{i \in I_x}\sum_{j \in J_x} \lambda_{ij}(x) = 1$ and
	\begin{equation}\label{**eq34}
		u(x)= t(x)+\beta \sum\limits_{j\in J_x}  \sum\limits_{i\in I_x}\lambda_{ij}(x) \nabla h_j(x,\xi_i).
	\end{equation}
If it is taken that $\lambda_{ij}(x)=0$ for $i\not\in I_x$ and $j\not\in J_x$, then by equation (\ref{**eq34}), it is obtain that
\begin{equation}\label{eq34}
u(x)= t(x)+\beta \sum\limits_{j\in \Lambda}  \sum\limits_{i\in \bar{\Lambda}}\lambda_{ij}(x) \nabla h_j(x,\xi_i)
\end{equation}
with  $\sum\limits_{j\in \Lambda}  \sum\limits_{i\in \bar{\Lambda}}\lambda_{ij}(x) =1.$
Combining (\ref{eq33}) and (\ref{eq34}), it follows that
$ \langle \beta \sum\limits_{j\in \Lambda}  \sum\limits_{i\in \bar{\Lambda}}\lambda_{ij}(x) \nabla h_j(x,\xi_i)+t(x), \; t-t(x) \rangle \geq0$ for each $t\in D-\{x\}.$ Since $t\in D-\{x\}$ and if it is considered that $z=t+x \in D,$ then the later expression can be written as \begin{equation}\label{3.5} \langle x- \beta \sum\limits_{j\in \Lambda}  \sum\limits_{i\in \bar{\Lambda}}\lambda_{ij}(x) \nabla h_j(x,\xi_i)-(x+t(x)), \; z-(x+t(x)) \rangle \leq0,
\end{equation} 
for each $z\in D.$ Therefore, by projection theorem \big(see Proposition 2.2.1 in \cite{bertsekas2003convex}\big), $x+t(x)$ is the projection of $x- \beta \sum\limits_{j\in \Lambda}  \sum\limits_{i\in \bar{\Lambda}}\lambda_{ij} \nabla h_j(x,\xi_i)$  on $D,$ then $$x+t(x)=\text{Proj}_D\bigg(x- \beta \sum\limits_{j\in \Lambda}  \sum\limits_{i\in \bar{\Lambda}}\lambda_{ij} \nabla h_j(x,\xi_i)\bigg),$$ 
where $\text{Proj}_D(y)=\underset{u\in D}{\mathrm{argmin}}\|u-y\|$ and $y=x- \beta \sum\limits_{j\in \Lambda}  \sum\limits_{i\in \bar{\Lambda}}\lambda_{ij} \nabla h_j(x,\xi_i).$ Since $D$ is closed and convex, then $\text{Proj}_D(y)$ is well defined. It can be conclude that there exists $\lambda\in\mathbb{R}^{m\times p}_{\geq}$ with $\lambda_{ij}\geq0,$ $\sum\limits_{j\in \Lambda}  \sum\limits_{i\in \bar{\Lambda}}\lambda_{ij} =1$ such that 
\begin{equation}\label{projectdescent}
t(x)=\text{Proj}_D\bigg(x- \beta \sum\limits_{j\in \Lambda}  \sum\limits_{i\in \bar{\Lambda}}\lambda_{ij} \nabla h_j(x,\xi_i)\bigg)-x.
\end{equation}
Thus, $t(x)$  is the solution for the subproblem $\text{S}_P(x)$ defined in \eqref{eq32}. Furthermore, the optimal value of the subproblem \eqref{eq32} is given by
\begin{equation}\label{theta}
\Omega(x)= \beta\Theta_{x}(t(x))+\frac{1}{2}\|t(x)\|^2.
\end{equation} Also, $t=\hat 0=(0,0,\ldots,0)\in D-\{x\},$ then it follows that 
\begin{align*}	
\Omega(x) &\leq \displaystyle \max_{j\in \Lambda}\{\max_{i\in \bar\Lambda}h_j(x,\xi_i) + \max_{i\in \bar\Lambda}\nabla h_j(x,\xi_i)^T\hat 0 \\
&+\max_{i\in \bar\Lambda} \frac{1}{2}\hat 0^TH_j(x,\xi_i)\hat 0-\max_{i\in \bar\Lambda}H_j(x)\}+\frac{1}{2}\|\hat 0\|^2.\end{align*} Hence, $\Omega(x)\leq0$.
\begin{remark} 
	If $p=1$, i.e., $U=\{\xi_1\}$, then it is obtained that $\Theta_{x}(t)= \displaystyle \max_{j\in\Lambda}\max_{i\in \bar{\Lambda}} \{ h_j(x,\xi_i)+\nabla h_j(x,\xi_i)^T t - H_j(x)\} = \displaystyle \max_{j\in\Lambda} \{ \nabla h_j(x,\xi_1)^T t\}$. In this case, subproblem (\ref{eq32}) reduces to the subproblem for PGM for deterministic multiobjective optimization problems. Hence, it can be said that subproblem (\ref{eq32}) is just an extension of the subproblem used for the projected gradient descent direction in the PGM for multiobjective optimization problems given by Drummond et al. \cite{drummond2004projected}.
	
\end{remark}
\begin{remark} 
If $p=1$ and $m=1,$ then $\Theta_{x}(t)= \nabla h_j(x,\xi_1)^Tt$. Thus the subproblem (\ref{eq32}) reduced to subproblem used for projected gradient descent direction in PGM for scalar optimization problem.
\end{remark}
\par In Theorem \ref{3.2t}, it is shown that $\Omega(x)$ is a continuous function on $D.$
\begin{thm}\label{3.2t}
Let at any $x\in D,$ $\Omega(x)$ be the optimal solution of the subproblem  $\displaystyle \min_{t\in {D-\{x\}}}\left(\beta \Theta_{x}(t)+\tfrac{1}{2}\|t\|^2 \right),$ then $\Omega(x)$ is continuous function on $D.$ 
\end{thm}
\begin{proof}
By equation (\ref{theta}), $\Omega(x)= \beta\Theta_{x}(t(x))+\frac{1}{2}\|t(x)\|^2,$ where $\Theta_{x}(t)= \displaystyle \max_{j\in\Lambda}\max_{i\in \bar{\Lambda}}  \{ h_j(x,\xi_i)+\nabla h_j(x,\xi_i)^Tt-H_j(x)\}$ is continuous and $\beta>0$ is a parameter. The continuity of $\Omega(x)$ will be demonstrated by employing the sequential criterion of continuity. For that let $x\in D$ and $\{x^k\}$	such that $x^k\in D$ and $x^k\to x$ as $k\to \infty.$ As $t(x)\in D-\{x\},$ it follows that $t(x)+x-x^k\in D-\{x^k\}.$ Now by definition of $\Omega$ at $x^k,$ it is obtained that
\begin{align*}
	\Omega(x^k)&\leq \beta\Theta_{x^k}(t(x)+x-x^k)+\frac{1}{2}\|t(x)+x-x^k\|^2\\
	&\leq \beta\Theta_{x^k}(t(x))+\beta\Theta_{x^k}(x-x^k) +\frac{1}{2}\|t(x)\|^2+\frac{1}{2}\|x-x^k\|^2+\langle t(x),x-x^k\rangle.
\end{align*}
After taking $\limsup\limits_{k\rightarrow\infty},$ it follows that 
\begin{align}
	\limsup\limits_{k\rightarrow\infty}	\Omega(x^k)&\leq \beta\Theta_{x}(t(x))+\frac{1}{2}\|t(x)\|^2= \Omega(x)\implies \limsup\limits_{k\rightarrow\infty}	\Omega(x^k)\leq  \Omega(x)\label{lms}. 
\end{align}
Since $t(x^k)\in D-\{x^k\},$ it is obtained that $t(x^k)+x^k-x\in D-\{x\}.$ Again by $\Omega$ at $x,$ it follows that
\begin{align*}
	\Omega(x)&\leq \beta\Theta_{x}(t(x^k)+x^k-x)+\frac{1}{2}\|t(x^k)+x^k-x\|^2\\
	&\leq \beta\Theta_{x}(t(x^k))+\beta\Theta_{x}(x^k-x) +\frac{1}{2}\|t(x^k)\|^2+\frac{1}{2}\|x^k-x\|^2+\langle t(x^k),x^k-x\rangle.
\end{align*}
By taking the $\liminf\limits_{k\rightarrow\infty},$ it is obtained that
\begin{align}
	\Omega(x)&\leq\liminf\limits_{k\rightarrow\infty}[\beta\Theta_{x}(t(x^k))+\frac{1}{2}\|t(x^k)\|^2]\nonumber\\
	&\leq\liminf\limits_{k\rightarrow\infty}[\beta\Theta_{x}(t(x^k))+\frac{1}{2}\|t(x^k)\|^2]\nonumber\\
	&\leq\liminf\limits_{k\rightarrow\infty}[\beta\Theta_{x}(t(x^k))+\beta\Theta_{x^k}(t(x^k))-\beta\Theta_{x^k}(t(x^k))+\frac{1}{2}\|t(x^k)\|^2]\nonumber\\
	&\leq\liminf\limits_{k\rightarrow\infty}[\Omega(x^k)+ \beta(\Theta_{x}(t(x^k))-\beta\Theta_{x^k}(t(x^k))]\nonumber\\
	&\leq\liminf\limits_{k\rightarrow\infty}\Omega(x^k)\label{liminf}.
\end{align}
Now by the inequalities (\ref{lms}) and (\ref{liminf}), it follows that $\Omega(x)\leq\liminf\limits_{k\rightarrow\infty}\Omega(x^k)\leq\limsup\limits_{k\rightarrow\infty}	\Omega(x^k)\leq \Omega(x).$ Consequently, \[\lim_{k\to\infty}\Omega(x^k)=\Omega(x)\] and hence $\Omega(x)$ is continuous.
\end{proof}
\par In the next theorem, denoted as Theorem \ref{Thm 3.2.}, the characterization of the Pareto critical point in terms of $t(x)$ and $\Omega(x)$ will be discussed. Additionally, some features of $t(x)$ and $\Omega(x)$ will be provided. It is also demonstrated that $t(x)$ is the descent direction for $H$ at $x$.
\begin{thm}\label{Thm 3.2.} Let $\Omega(x)$ and $t(x)$ are the optimal value and optimal solution of the descent direction finding subproblem $(\ref{eq32})$, respectively, which are given in the equations $(\ref{eq31})$ and $(\ref{theta})$. Then the following results will hold:
\begin{enumerate}
	\item $ t(x)$ is bounded on any compact subset $C$ of  $D.$
	\item The next four conditions are equivalent:
	\begin{enumerate}
		\item The point $x$ is not a Pareto critical point.
		\item $\Omega(x)<0.$
		\item $t(x)\not =0.$
		\item$t(x)$ is a descent direction for $H$ at $x$.
	\end{enumerate}
\end{enumerate}
In particular, $x$ is a  Pareto critical point $\Longleftrightarrow$ $\Omega(x)=0$.
\end{thm}
\begin{proof} 1. Let $C$ be a compact subset of $\mathbb{R}^n$. Since $h_j(x,\xi_i)$ is continuously differentiable for each $j\in \Lambda$ and $\xi_i\in U, ~i\in\bar{\Lambda}$, then $h_j(x,\xi_i)$ is bounded on every compact set. So, for all $x\in C$, $\xi_i\in U,~i\in\bar{\Lambda}$ and $j \in \Lambda,$ by equation $(\ref{projectdescent})$  $t(x)$ is bounded on the compact set $C$.\\ \\
$2.~(a)\implies(b):$ \\
Assume that $x$ is not a Pareto critical point. Then there exists $\bar t$ such that\begin{center}
	$\nabla h_j(x,\xi_i)^T \bar t(x) < 0, ~\forall ~j \in \Lambda$ and $i \in I_j(x)$.
\end{center}
Since $\Omega(x)$ is the optimal value for the subproblem \eqref{eq32}, then for all $\delta >0,$ it follows that
\begin{align*}
	\Omega(x) &\leq \beta \displaystyle \max_{j\in \Lambda}\max_{i\in \bar{\Lambda}} \{f_j(x,\xi)+\nabla h_j(x,\xi_i)^T(\delta \bar t)-H_j(x)\}+\frac{1}{2}\|\delta \bar t\|^2\nonumber\\
	&\leq \beta  \displaystyle \max_{j\in \Lambda}\{\max_{i\in \bar{\Lambda}} h_j(x,\xi_i)+ \max_{i\in \bar{\Lambda}}\nabla h_j(x,\xi_i)^T (\delta\bar t)-\max_{i\in \bar{\Lambda}} H_j(x)\}+\frac{1}{2}\| \delta\bar t\|^2 \nonumber\\
	&=\beta \displaystyle \max_{j\in \Lambda}\max_{i\in \bar{\Lambda}}\nabla h_j(x,\xi_i)^T (\delta\bar t)+\frac{1}{2}\| \delta\bar t\|^2 \\
	&=\beta \delta\big(\displaystyle \max_{j\in \Lambda}\max_{i\in \bar{\Lambda}}\nabla h_j(x,\xi_i)^T (\bar t)+\frac{1}{2}\delta\| \bar t\|^2 \big).
\end{align*}
For small enough $\delta>0$, the right-hand side of the above inequality will be negative because of $\nabla h_j(x,\xi_i)^Tt(x) < 0$ for all $j\in\Lambda$ and $i\in \bar{\Lambda}$.  Thus, $\Omega(x)
< 0$.\\ 
$(b) \implies (c):$\\  Given that $\Omega(x)$ represents the optimal value of the subproblem \eqref{eq32} and is negative according to condition $(b)$, it follows that $t(x) \neq 0$. Because if $t(x)  = 0,$ then $\Omega(x)$ will be zero which is not possible from $(b)$. Hence, if $\Omega(x) < 0,$ then $t(x) \not = 0$.\\
$(c) \implies (d):$ \\ Let $t(x) \not = 0$. Then $\Omega (x) \not =0$. Since $\Omega(x) \leq 0$, then $\Omega (x)<0$. Thus,
\allowdisplaybreaks
\begin{align*}
	\Omega(x)& =\beta \displaystyle \max_{j\in \Lambda}\max_{i\in\bar{\Lambda}}\{h_j(x,\xi_i)+\nabla h_j(x,\xi_i)^Tt(x)-H_{j}(x)\} +\frac{1}{2}\|t(x)\|^2 < 0.
\end{align*}Which implies 
\begin{align*}
	&\beta \displaystyle \max_{j\in \Lambda}\{\max_{i\in\bar{\Lambda}} h_j(x,\xi_i)+\max_{i\in \bar{\Lambda}}\nabla h_j(x,\xi_i)^Tt(x)-\max_{i\in \bar\Lambda}H_j(x)\}+\frac{1}{2}\|t(x)\|^2 < 0\\
	&\implies \beta\max_{j\in \Lambda}\max_{i\in \bar{\Lambda}}\nabla h_j(x,\xi_i)^Tt(x)< 0.	
\end{align*} \text{As $\beta>0$,} then
\begin{align*}
	&\nabla h_j(x,\xi_i)^Tt(x)< 0, ~\forall ~ j\in \Lambda ~\text{and}~ i\in \bar{\Lambda}\\
	&\implies\nabla h_j(x,\xi_i)^Tt(x)< 0, ~\forall ~ j\in \Lambda ~\text{and}~ i\in I_j(x)\\
	&\implies t(x)~ \text{is a descent direction for}~ H~ \text{at} ~x.
\end{align*}
$(d) \implies (a)$:  \\ Since $t(x)$ is a descent direction for $H$ at $x$, it follows that \begin{equation*}
	\nabla h_j(x,\xi_i)^Tt(x)< 0, ~\forall ~ j\in \Lambda ~\text{and}~ i\in I_j(x).
\end{equation*}
Then $x$ is not a Pareto critical point.\\
Also,  if $\Omega(x)<0$, then $t(x) \not= 0$. Moreover, if $t(x) \not= 0$, then $\Omega(x)<0$. Thus, $x$ is Pareto critical point $\Longleftrightarrow$ $\Omega(x) =0$.
\end{proof}
In subsection \ref{ss3.1}, the descent direction for the projected gradient method for $OWC_{P(U)}$ has been determined. To proceed with writing the projected gradient algorithm for $OWC_{P(U)}$, the next requirement is to determine the step length size.
\subsection{\textbf{Armijo-type inexact line search for the step length size}}\label{ss3.2}
To determine an appropriate step length that guarantees a  significant reduction in each objective function, an inexact line search technique is developed in this subsection. For the same purpose, an auxiliary function $H^*_j(x,t)$ is considered, defined as 
\begin{equation}\label{8}
H^*_j(x,t)= \max_{i\in\bar\Lambda}\{ \beta h_j(x,\xi_i) + \beta \nabla h_j(x,\xi_i)^Tt \}-\beta H_j(x),~ j\in \Lambda.~~~
\end{equation}
Note that equation (\ref{8}) implies
\begin{equation}\label{8**}
H^*_j(x,t) \leq \beta\displaystyle\max_{i\in \bar{\Lambda}} h_j(x,\xi_i) + \beta\displaystyle\max_{i\in \bar{\Lambda}}\nabla h_j(x,\xi_i)^Tt -\beta H_j(x),~ j\in \Lambda,~~
\end{equation}
Since $\Omega(x)\leq0$, it follows that
\begin{center}
$\Omega(x)	= \beta\displaystyle \max_{j\in\Lambda}\max_{i\in \bar{\Lambda}}  \{h_j(x,\xi_i)+\nabla h_j(x,\xi_i)^Tt-H_j(x)\} +\frac{1}{2} \|t\|^2\leq0,$
\end{center}
which gives
\begin{equation}\label{3.9ls}
\beta\displaystyle \max_{j\in\Lambda}\max_{i\in \bar{\Lambda}}  \{h_j(x,\xi_i)+\nabla h_j(x,\xi_i)^Tt-H_j(x)\}\leq-\frac{1}{2} \|t\|^2.
\end{equation}
Now from inequality (\ref{3.9ls}), it follows that
\begin{equation*}
\beta(h_j(x,\xi_i)+\nabla h_j(x,\xi_i)^Tv -H_{j}(x))\leq 0,~ ~\text{for each}  ~i ~\text{and} ~j, 
\end{equation*}
which implies 
\begin{equation}\label{9}
\beta\displaystyle\max_{i\in \bar{\Lambda}} h_j(x,\xi_i) + \beta\displaystyle\max_{i\in \bar{\Lambda}}\nabla h_j(x,\xi_i)^Tt -\beta H_j(x)\leq0, ~\text{for each} ~j.
\end{equation}
After combining the inequalities $ (\ref{8})$, $(\ref{8**})$ and $(\ref{9}),$ it is obtained that
\begin{equation}\label{11***}
H^*_j(x,t) \leq 0.
\end{equation}
Since $h_j(x,\xi_i)$ is convex and continuously differentiable function and hence it is an upper uniformly differentiable (see Definition 2.1, p.159 in \cite{bazaraa1982algorithm}). Since $h_j(x, \xi_i)$ is upper uniformly differentiable function, then there exists $k_i$ such that
\begin{center}
$	h_j(x+t,\xi_i) \leq h_j(x,\xi_i)+\nabla h_j(x,\xi_i)^Tt + \frac{1}{2} k^j_i \|t\|^2,~~\text{for each},~ i\in \bar\Lambda ~\text{and} ~j\in\Lambda.$
\end{center}
Also,
\begin{align*}
h_j(x+t,\xi_i) &\leq h_j(x,\xi_i)+\nabla h_j(x,\xi_i)^Tt + \frac{1}{2} k_i \|t\|^2 \\ &\leq \beta\max_{i\in \bar\Lambda} \{ h_j(x,\xi_i)+\nabla h_j(x,\xi_i)^Tt \}+ \frac{1}{2} K^j \|t\|^2,
\end{align*}
where $K^j= \displaystyle \max_{i\in \bar\Lambda}{k^j_i}$ and $\beta\geq1$. Now, from the equation (\ref{8}), it is established that
\begin{equation*}
h_j(x+t,\xi_i) \leq 	H^*_j(x,t) + \beta H_j(x)+\frac{1}{2} K^j \|t\|^2,
\end{equation*}
which holds for each $i\in \bar{\Lambda}$, and $j\in\Lambda$. Therefore,
\begin{align}
&\max_{i\in  \bar{\Lambda}}h_j(x+t,\xi_i) \leq 	H^*_j(x,t) + \beta H_j(x)+\frac{1}{2} K^j \|t\|^2 \notag 
\\ 
\text{i.e.,} ~&H_j(x+t)\leq H^*_j(x,t) + \beta H_j(x)+\frac{1}{2}K^j \|t\|^2\label{10}.
\end{align}
The following  two results related to continuous approximation are true for all $j\in\Lambda$ \cite{bazaraa1982algorithm}:
\begin{enumerate}[($i$)]
\item\label{m} 	$H^*_{j}(x,\lambda t)\leq \lambda H^*_{j}(x,t),~\forall \lambda \in [0,1]$.
\item $H_j(x+\lambda t)\leq  H_j(x) +\lambda H^*_j(x,t)+\frac{1}{2}\lambda^2 K^j \|t\|^2.$
\end{enumerate}
Now put $t = \alpha t$ in equation (\ref{10}), then with the help of $(\ref{m})$, it is obtained that
\begin{equation}\label{11}
H_j(x+\alpha t)\leq \alpha H^*_j(x,t) + \beta H_j(x)+\frac{1}{2}\alpha^2K^j \|t\|^2,
\end{equation}
where $\alpha>0$ is sufficiently small.\\
Now, if $t\not=0$ then by inequality (\ref{11***}), it follows that  
\begin{center}
$H^*_j(x,t)< 0. $
\end{center}
For any $\eta \in (0,1),$ it is implies that
\begin{equation}\label{12}
H^*_j(x,t) < \eta  H^*_j(x,t).
\end{equation}
As $1>\alpha>0$ and $\beta\geq1$, if $\beta=1$, then by inequalities (\ref{11}) and (\ref{12}), it follows that 
\begin{center}$H_j(x+\alpha t)\leq  H_j(x)+  \alpha \eta H^*_j(x,t)+\frac{1}{2}\alpha^2K^j \|t\|^2.$
\end{center}
Since $\alpha$ is sufficiently small and positive, then
\begin{equation}\label{13}
H_j(x+\alpha t)\leq  H_j(x)+  \alpha \eta H^*_j(x,t).
\end{equation}
Inequality (\ref{13}) represents the step size rule for the projected gradient descent algorithm for the $OWC_{P(U)}$.
\par Having determined a descent direction for PGM and an inexact line search approach to find the step length size, the next step involves writing the projected gradient descent algorithm for $OWC_{P(U)}$ as follows:\\ \ \\ \ \\ 
\subsection{\textbf{Algorithm \big(Projected gradient descent algorithm for $OWC_{P(U)}$\big)}}\label{algo1}
\begin{enumerate}[{Step} 1]
\item Choose $\epsilon>0,~\beta\geq1,$~$\eta \in (0,1),$~and~$x^0\in D$. Set $k :=0.$
\item\label{st2} Solve $\text{S}_P(x^k)$ (defined in equation \ref{eq32}) and find $t^k$~and~$\Omega({x^k})$.
\item\label{st3} If  $|\Omega (x^k)|<\epsilon$, then stop. Otherwise proceed to Step \ref{st4}.
\item \label{st4} Choose $\alpha_k$ as the largest $\alpha \in \{ \frac{1}{2^r} : r=1,2,3,\ldots\} $ satisfying inequality (\ref{13})
such that
\begin{equation}\label{amj1}
	H_j(x^k+\alpha t^k)\leq H_j(x^k)+\alpha \eta H^*_j(x^k,t^k),~~\forall~j\in \Lambda, ~~~
\end{equation}
where  $H^*_j(x^k,t^k) = \displaystyle \max_{i \in \bar \Lambda} \{h_j(x^k,\xi_i) + \nabla h_j(x^k,\xi_i)^Tt^k\}- H_j(x^k)$.
\item\label{step5} Define $x^{k+1}:= x^k + \alpha_{k} t^k$, update $k:=k+1$ and go to Step 2.
\end{enumerate}
\textbf{Comments related to Algorithm $\ref{algo1}$ :}  
The well-definedness of Algorithm \ref{algo1} depends on Step~\ref{st2}, Step~\ref{st3}, Step~\ref{st4}, and Step~\ref{step5}. As the objective function $\Theta_{x^k}(t)+\frac{1}{2}\|t\|^2$ of the subproblem $\text{S}_P(x^k)$ is strongly convex, it has a unique minimizer, namely $t(x^k)=t^k$. The optimal value $\Omega(x^k)$ of the subproblem $\text{S}_P(x^k)$ can then be calculated, ensuring that Step~\ref{st2} is well defined. In Step~\ref{st3}, $|\Omega(x^k)|$ can be easily calculated, and $|\Omega(x^k)|<\epsilon$ is considered the stopping criterion for the algorithm. In Step~\ref{st4}, the step length size $\alpha_k$ at $x^k$ can be considered the largest $\alpha\in \left\{\frac{1}{2^r} : r=1,2,3,\ldots\right\}$ such that (\ref{amj1}) is satisfied. Using $\alpha_k$ and $t^k$, the new iterative scheme $x^{k+1}=x^k+\alpha_k t^k$ can be written. By Lemma \ref{l3.2}, it can be observed that $\{x^k\}\subset D$, ensuring that Step~\ref{step5} is well defined.
\par After presenting the Algorithm \ref{algo1} in textual form, a visual representation of the Algorithm \ref{algo1} is presented in the form of a flowchart. This flowchart illustrates the sequential steps outlined in the Algorithm \ref{algo1} and offers a more intuitive depiction of the Algorithm \ref{algo1}'s process.
\subsubsection{\textbf{Flowchart of Algorithm \ref{algo1}}}
\tikzstyle{startstop} = [rectangle, rounded corners, minimum width=3cm, minimum height=1cm,text centered, draw=black, fill=red!30]
\tikzstyle{process} = [rectangle, minimum width=3cm, minimum height=1cm, text centered, draw=black, fill=orange!30]
\tikzstyle{decision} = [diamond, minimum width=3cm, minimum height=1cm, text centered, draw=black, fill=green!30]
\tikzstyle{arrow} = [thick,->,>=stealth]

\begin{tikzpicture}[node distance=2cm]
	
	\node (start) [startstop] {Start};
	\node (init) [process, below of=start] {Initialization: Choose $\epsilon>0,~\beta\geq1,$~$\eta \in (0,1),$~and~$x^0\in D$. Set $k :=0.$};
	\node (solve) [process, below of=init] {Solve $\text{S}_P(x^k)$ for $t^k$~and~$\Omega({x^k})$};
	\node (check) [decision, below of=solve, yshift=-0.7cm] {If  $|\Omega (x^k)|<\epsilon$ or $\|t^k\|<\epsilon$};
	\node (stop) [startstop, right of=check, xshift=3cm] {Stop};
	\node (choose) [process, below of=check, yshift=-1.0cm] {Choose $\alpha_k$ as the largest $\alpha \in \{ \frac{1}{2^r} : r=1,2,3,\ldots\} $ satisfying inequality (\ref{13}).};
	\node (update) [process, below of=choose] {Update Solution: $x^{k+1}:= x^k + \alpha_{k} t^k$};
	\node (increment) [process, below of=update] {Increment Iteration: $k:=k+1$};
	\node (loop) [startstop, below of=increment] {Loop Back to Solve ${S_P}(x^k)$};
	
	\draw [arrow] (start) -- (init);
	\draw [arrow] (init) -- (solve);
	\draw [arrow] (solve) -- (check);
	\draw [arrow] (check) -- node[anchor=south] {Yes} (stop);
	\draw [arrow] (check) -- node[anchor=east] {No} (choose);
	\draw [arrow] (choose) -- (update);
	\draw [arrow] (update) -- (increment);
	\draw [arrow] (increment) -- (loop);
	\draw [arrow] (loop.west) -- ++(-4,0) |- (solve.west);
\end{tikzpicture}
\subsection{\textbf{Convergence analysis of Algorithm \ref{algo1}}}\label{ss3.3}
It is obvious that if Algorithm \ref{algo1} has a finite number of iterations, then the last iterative point is approximately a Pareto critical point, and therefore it is a Pareto optimal solution for $H$. Consequently, when Algorithm \ref{algo1} generates an infinite sequence, it is essential to consider the convergence analysis. In view of this, it is presumed that $\{x^k\}$, $\{t^k\}$, and $\{\alpha_k\}$ are infinite sequences generated by Algorithm \ref{algo1} for $OWC_{P(U)}$. To demonstrate that any accumulation point of $\{x^k\}$ is a Pareto critical point for $H$, it is first shown that every term of the sequence generated by Algorithm \ref{algo1} lies in $D$.

\begin{lem}\label{l3.2}
Assuming that the sequence $\{x^k\}\subset \mathbb{R}^n$ is created by Algorithm \ref{algo1}, then  $x^k\in D$ for all $k.$
\end{lem}
\begin{proof}
The proof of this lemma can be provided by induction. Starting with $k=0$, it is known that $x^0\in D$. Now, assuming $x^k\in D$, it needs to be demonstrated that $x^{k+1}\in D$. From Step \ref{st4}, it is observed that $\alpha\in(0,1]$. Since $t^k=t(x^k)\in D-{x^k}$, then for some $p^k\in D$ there exists $t^k=p^k-x^k$. Then, by the convexity assumption of $D$, it is clear that $x^{k+1}=x^k+\alpha_kt^k=x^k+\alpha_k(p^k-x^k)=(1-\alpha_k)x^k+\alpha_kp^k\in D$. Therefore, by induction, $x^k\in D$ for all $k$.
\end{proof}
\begin{thm}\label{th3.3}
Assume that the sequence $\{x^k\}$ is created by the Algorithm \ref{algo1}. Any point which is an accumulation point of $\{x^k\}$ is a Pareto critical point for $H$.
\end{thm}
\begin{proof}
Let  $x^*$ be an accumulation point  of the sequence $\{x^k\}$. Consider $t(x^*)$ and $\Omega(x^*)$ that are given by $t(x^*)	= \underset{t\in \mathbb{R}^n}{\mathrm{argmin}} \left\{  \beta\displaystyle \max_{j\in\Lambda}\max_{i\in \bar{\Lambda}}  \{ h_j(x^*,\xi_i)+\nabla h_j(x^*,\xi_i)^Tt-H_j(x^*)\} +\frac{1}{2} \|t\|^2\right\}$
\text{and}~~$\Omega(x^*)=  \beta \displaystyle \max_{j\in\Lambda}\max_{i\in \bar{\Lambda}}  \{ h_j(x^*,\xi_i)+\nabla h_j(x^*,\xi_i)^Tt(x^*)-H_j(x^*)\} +\frac{1}{2} \|t(x^*)\|^2$
respectively, i.e., the solution and optimal value of the subproblem \eqref{eq32} at $x=x^*$.
It is known that $H(x^k)$ is $\mathbb{R}^m$-decreasing (i.e., component-wise decreasing). So, it is obtained that
\begin{center}
	$\lim\limits_{k\rightarrow\infty}H(x^k) = H(x^*).$
\end{center}
Therefore,
\begin{center}
	$\lim\limits_{k\rightarrow\infty}\|H(x^k) - H(x^*)\|=0.$
\end{center}
Component-wise it is obtained that
\begin{center}
	$ \lim\limits_{k\rightarrow\infty}H_{j}(x^k) = H_{j}(x^*), ~\text{for all}~~j=1,2,\ldots,m.$
\end{center}
We know from Step~\ref{st4} that
\begin{center}
	~~~~~~$	H_j(x^k+\alpha_k t^k)\leq H_j(x^k)+\alpha_k \eta H^*_j(x^k,t^k)$\\
	\text{and}~~~~$	H_j(x^k)-H_j(x^k+\alpha_k t^k) \geq -\alpha_k \eta H^*_j(x^k,t^k).$
\end{center}
Therefore,
\begin{equation}\label{14}
	\displaystyle \lim\limits_{k\rightarrow\infty}\alpha_k H^*_{j}(x^k,t^k)=0, ~\text{for all}~ j \in\Lambda.
\end{equation}
Since $\alpha_{k}\in(0,1],~\forall ~k$, the following two possibilities arise:
\begin{align}\label{15}
	& \limsup\limits_{k\rightarrow\infty} \alpha_{k}>0 ~ \text{ or } \\
	\label{16}
	& \limsup\limits_{k\rightarrow\infty} \alpha_{k}=0.
\end{align}
Assume that (\ref{15}) holds, a convergent subsequence $\{x^{k_l}\}$ of the sequence $\{x^k\}$ exists that converges to $x^*,$ and $\alpha^*$ such that $\lim\limits_{l\rightarrow\infty}\alpha_{k_l}=\alpha^*$. So, using (\ref{14})
it is established that
\begin{align*}
	0&=\lim\limits_{l\rightarrow\infty}\alpha_{k_l}H^*_j(x^{k_l},t^{k_l})\\
	&\leq \lim\limits_{l\rightarrow\infty} \bigg( \beta(h_j(x,\xi_i)+\displaystyle \max_{j\in\Lambda}\max_{i\in J_j(x^{k_l})}  \nabla h_j(x^{k_l},\xi_i)^Tt^{k_l}-H_j(x))+\tfrac{1}{2} \|t(x^{k_l})\|^2\bigg),\\
	&=\Omega(x^{k_l})
\end{align*}
and therefore $0\leq \lim\limits_{l\rightarrow\infty}\Omega(x^{k_l}).$
Since $\Omega(x)\leq 0,$ for each $x$, it is obtained that
$0\leq\Omega(x^*)\leq 0,$ which implies $\Omega(x^*)= 0$ and hence $x^*$ is a Pareto critical point.\\
Now, it is assumed that (\ref{16}) holds. By  Theorem \ref{Thm 3.2.}, $t$ is bounded on any compact set, and $t^k=t(x^k), ~\forall~ k$ and $i\in I_j(x)$. It follows that $\{x^k\}$ has a convergent subsequence.~It also follows that $\{t^k\}$ has a bounded subsequence. Therefore, there exists $t^* \in \mathbb{R}^n$ and a subsequence $\{t^{k_l}\}_z$ such that $\lim\limits_{l\rightarrow\infty}t^{k_l}=t^*$ and $\lim\limits_{l\rightarrow\infty}\alpha_{k_l}=0.$~Note that
\begin{center}
	$\beta\displaystyle \max_{j\in\Lambda}\max_{i\in\bar{\Lambda}}  \{h_j(x,\xi_i)+\nabla h_j(x^{k_l},\xi_i)^Tt^{k_l}-H_j(x)\}\leq \Omega(x^{k_l})<0, ~\forall ~l.$
\end{center}
So, as $l\to\infty$ it follows that
\begin{equation}\label{17*}
	\beta	\max_{j\in\Lambda}\max_{i\in \bar{\Lambda}} \{h_j(x^*,\xi_i)+ \nabla h_j(x^*,\xi_i)^Tt^*-H_j(x^*)\}\leq \Omega(x^*)\leq0.
\end{equation}
Since $\alpha_{k_l} \rightarrow 0$ for $l$ large enough, by Archimedian  property, $\alpha_{k_l}<\frac{1}{2^r}$, which means that the Armijo-like line search rule is not satisfied for $\alpha = \frac{1}{2^r},$ i.e.,
\begin{center}
	$	H_j(x^{k_l}+\frac{1}{2^r} t^{k_l})\not\leq H_j(x^{k_l})+ \frac{1}{2^r}\eta H^*_j(x^{k_l},t^{k_l}).$
\end{center}
So for all $j$, there exists $j=j(k_l)\in \Lambda$ such that
\begin{center}
	$	H_j(x^{k_l}+\frac{1}{2^r} t^{k_l})\geq H_j(x^{k_l})+ \frac{1}{2^r}\eta H^*_j(x^{k_l},t^{k_l}).$
\end{center}
Since  $\{j(k_l)\}_l \in \Lambda$, there exist a subsequence $\{ k_{l_z}\}_z$
and an index $j_0$ such that\\ $j_0=j(k_{l_z}),~\forall~ z = 1, 2, 3,\ldots$  and
\begin{center}
	$	H_j(x^{k_{l_z}}+\frac{1}{2^r} t^{k_{l_z}}) \geq H_j(x^{k_{l_z}}) + \frac{1}{2^r}\beta H^*_j({} x^{k_{l_z}},t^{k_{l_z}}). $
\end{center}
Taking limit $z \rightarrow \infty$ in the above inequality, it is obtained that
\begin{center}
	$H_j(x^*+\frac{1}{2^r} t^*) \geq H_j(x^*) + \frac{1}{2^r}\eta H^*_j( x^*,t^*).$
\end{center}
Given that this inequality is applicable to all positive integers $r$ and for $j_0$ (depends on $r$), then from (\ref{13})
\begin{center}
	$	H_j(x+\alpha t)\leq  H_j(x)+  \alpha \eta H^*_j(x,t).$
\end{center}
It follows that
\begin{center}
	$H^*_j( x^*,t^*)\not <0$.
\end{center}
So,
\begin{equation}\label{18*}
	\beta	\max_{j\in\Lambda}\max_{i\in \bar{\Lambda}} \{h_j(x^*,\xi_i)+ \nabla h_j(x^*,\xi_i)^Tt^*-H_j(x^*)\}\geq0.
\end{equation}
From inequalities (\ref{17*}) and (\ref{18*}), $\Omega (x^*)=0$. Consequently, it can be concluded that $x^*$ is a Pareto critical point.
\end{proof}
\subsubsection{\textbf{Global convergence of Algorithm \ref{algo1} (projected gradient algorithm (PGA))}}\label{s341} 
To demonstrate the global convergence of Algorithm \ref{algo1}, first, Lemma \ref{lem3.3} is proved.
\begin{lem}\label{lem3.3}
Let $H$ be a $\mathbb{R}^m_{\geq}$-convex function and $\sum\limits_{j\in \Lambda}  \sum\limits_{i\in \bar{\Lambda}} \nabla h_j(x,\xi_i)^T (x-y)\leq H(x)-H(y)$ for all $x\in D.$ Let $\{x^k\}$ be the sequence produced by Algorithm \ref{algo1}. Now, if for $x'\in D$ and $k\geq0,$ $H(x')\leq H(x^k),$ then $$\|x^{k+1}-x'\|^2\leq\|x^k-x'\|^2+2\beta \alpha_k|\langle \lambda^k, \; \sum\limits_{j\in \Lambda}  \sum\limits_{i\in \bar{\Lambda}} \nabla h_j(x,\xi_i)^Tt^k  \rangle |,$$ where $\lambda^k=\lambda(x^k)\in\mathbb{R}^{m\times p}_{\geq}$ is such that $\lambda_{ij}\geq0,$ $\sum\limits_{j\in \Lambda}  \sum\limits_{i\in \bar{\Lambda}}\lambda_{ij} =1$ and (\ref{projectdescent}) holds.
\end{lem}
\begin{proof}
By the Step~\ref{step5} of Algorithm \ref{algo1}, $x^{k+1}=x^k+\alpha_kt^k,$ it follows that
$x^{k+1}-x'=x^k+\alpha_kt^k-x',$ which implies 	\begin{equation}\label{fe}
	\|x^{k+1}-x'\|^2=\|x^k-x'\|^2+\alpha^2_k\|t^k\|^2-2\alpha_{k}\langle t^k,\;x'-x^k\rangle.~~~
\end{equation}
Also, it is given that $t^k=t(x^k)$, $\lambda^k=\lambda(x^k).$ Then by (\ref{3.5}) with $x=x^k,$ it is established that 
$ \langle \beta \sum\limits_{j\in \Lambda}  \sum\limits_{i\in \bar{\Lambda}}\lambda^k_{ij} \nabla h_j(x^k,\xi_i)+t^k, \; t-t^k \rangle \geq0$ for each $t^k\in D-\{x^k\}.$ Now take $t=x'-x^k\in D-\{x^k\}.$ Then by previous inequality, it is obtained that 
\begin{align*} \langle \beta \sum\limits_{j\in \Lambda}  \sum\limits_{i\in \bar{\Lambda}}\lambda^k_{ij} \nabla h_j(x^k,\xi_i), \; t \rangle+ \langle \beta \sum\limits_{j\in \Lambda}  \sum\limits_{i\in \bar{\Lambda}}\lambda^k_{ij} \nabla h_j(x^k,\xi_i), \; -t^k \rangle+\langle t^k,\;t-t^k\rangle\geq0.
\end{align*} Also,
\begin{align}\label{3.24} \langle \beta \sum\limits_{j\in \Lambda}  \sum\limits_{i\in \bar{\Lambda}}\lambda^k_{ij} \nabla h_j(x^k,\xi_i), \; t \rangle- \langle \beta \sum\limits_{j\in \Lambda}  \sum\limits_{i\in \bar{\Lambda}}\lambda^k_{ij} \nabla h_j(x^k,\xi_i), \; t^k \rangle -\|t^k\|^2\geq-\langle t^k,\;x'-x^k\rangle.
\end{align}
Inequality (\ref{3.24}) can be written as:
\begin{align}\label{3.25} \beta\langle  \lambda^k, \;\sum\limits_{j\in \Lambda}  \sum\limits_{i\in \bar{\Lambda}} \nabla h_j(x^k,\xi_i)^T (x'-x^k) \rangle- \beta\langle \lambda^k, \; \sum\limits_{j\in \Lambda}  \sum\limits_{i\in \bar{\Lambda}} \nabla h_j(x^k,\xi_i)^Tt^k \rangle& -\|t^k\|^2\geq-\langle t^k,\;x'-x^k\rangle.
\end{align}
Since by assumption $H$ is $\mathbb{R}^m_\geq$-convex, $\sum\limits_{j\in \Lambda}  \sum\limits_{i\in \bar{\Lambda}} \nabla h_j(x^k,\xi_i)^T (x'-x^k)\leq H(x')-H(x^k),$ for $x^k,x' \in D$  and $H(x')\leq H(x^k),$ then it follows that  $\sum\limits_{j\in \Lambda}  \sum\limits_{i\in \bar{\Lambda}} \nabla h_j(x^k,\xi_i)^T (x'-x^k)\leq0.$ Therefore,
\begin{equation}\label{3.26}
	\langle  \lambda^k, \;\sum\limits_{j\in \Lambda}  \sum\limits_{i\in \bar{\Lambda}} \nabla h_j(x^k,\xi_i)^T(x'-x^k) \rangle\leq \langle \lambda^k, \; H(x')-H(x^k) \rangle\leq0.~~
\end{equation} As $x^k$ is not a Pareto critical point, then $\nabla h_j(x^k,\xi_i)^Tt^k<0$ which follows that $\langle \lambda^k,\; \sum\limits_{j\in \Lambda}  \sum\limits_{i\in \bar{\Lambda}} \nabla h_j(x^k,\xi_i)^Tt^k\rangle<0.$ Since $\beta>0$ then by (\ref{3.25}) and (\ref{3.26}), it follows that 
\begin{equation}\label{3.27}
	-\langle t^k,\;x'-x^k\rangle\leq\beta|\langle  \lambda^k, \;\sum\limits_{j\in \Lambda}  \sum\limits_{i\in \bar{\Lambda}} \nabla h_j(x^k,\xi_i)^Tt^k \rangle|-\|t^k\|^2.
\end{equation}
Now by (\ref{fe}) and (\ref{3.27}), it is obtained that
\begin{align*}\|x^{k+1}-x'\|^2&\leq\|x^k-x'\|^2+\alpha^2_k\|t^k\|^2+2\beta\alpha_{k}|\langle  \lambda^k, \;\sum\limits_{j\in \Lambda}  \sum\limits_{i\in \bar{\Lambda}} \nabla h_j(x^k,\xi_i)^Tt^k \rangle|-2\alpha_k\|t^k\|^2.\end{align*}
As $\alpha_l\geq0$,  then by the above inequality result follows.
\end{proof}
The central idea for the convergence analysis is defined before proceeding to prove the next theorem. Remember that a sequence $\{y^k\},$ with $y^k\in \mathbb{R}^n$ for all $k$ is said be quasi-Fejer convergent to a set $Y\subset\mathbb{R}^n$ if for every $y\in Y$ there exists a sequence $\{\gamma_k\}$ with $\gamma_k\in\mathbb{R},$ $\gamma_k\geq0$ for all $k,$ and such that $\|y^{k+1}-y\|^2\leq\|y^k-y\|^2+\gamma_k$ for each $k=0,1,2,\ldots,$ with $\sum\limits^{\infty}_{k=0}\gamma_k<\infty.$ Also a well-known consequence associated with the quasi-Fejer convergent is that: If a $\{y^k\}$ is a quasi-Fejer convergent to a set $Y\subset\mathbb{R}^n$ then $\{y^k\}$ is bounded. Moreover, if $y\in Y$ is a cluster point of $\{y^k\},$ then $\displaystyle\lim_{k\to\infty}y^k=y.$
\begin{thm} Assume that $H$ is convex in component-wise sense and every $\mathbb{R}^m_\geq-$decreasing sequence $\{u^k\}$ in $H(D)$ is $\mathbb{R}^m_\geq-$ bounded below by any point of $H(D).$ Then the weak Pareto optimal solution is reached for any sequence $\{x^k\}$ generated using Algorithm \ref{algo1}.
\end{thm}
\begin{proof}
First, a set $S={x\in D:H(x)\leq H(x^k), k=0,1,2,\ldots,}$ is defined. Then, an element $x'\in S$ satisfying the assumption of this theorem is considered. It is also noted that $H$ is $\mathbb{R}^m_\geq$-convex; therefore, by Lemma \ref{lem3.3}, it follows that
\begin{equation}\label{3.29**}
	\|x^{k+1}-x'\|^2\leq\|x^k-x'\|^2+2\beta \alpha_k|\langle \lambda^k, \; \sum\limits_{j\in \Lambda}  \sum\limits_{i\in \bar{\Lambda}} \nabla h_j(x,\xi_i)^Tt^k  \rangle |.
\end{equation}
Since $\lambda^k_{ij}\geq0,$ $\sum\limits_{j\in \Lambda}  \sum\limits_{i\in \bar{\Lambda}}\lambda^k_{ij} =1$ then $0\leq\lambda^k_{ij}\leq1$ for all $i,j,$ and $k,$ also each $k,$ it can be written as 
\begin{equation}\label{delta}
	\lambda^k=\sum\limits_{j\in \Lambda}  \sum\limits_{i\in \bar{\Lambda}}\lambda^k_{ij}\delta_{ij},~ \text{where}~ \delta_{ij} =
	\begin{cases}
		0, &         \text{if }i\neq j,\\
		1, &         \text{if } i= j.
	\end{cases}
\end{equation}
By (\ref{3.29**}) and (\ref{delta}), it is obtained that
\begin{equation}\label{3.29}
	\|x^{k+1}-x'\|^2\leq\|x^k-x'\|^2+2\beta \alpha_k|\langle \delta_{ij}, \; \sum\limits_{j\in \Lambda}  \sum\limits_{i\in \bar{\Lambda}} \nabla h_j(x,\xi_i)^Tt^k  \rangle |.
\end{equation} 
Quasi-Fejer convergence of ${x^k}$ to the set $S$ is now demonstrated. If $\gamma_k=2\beta \alpha_k|\langle \delta_{ij}, \; \sum\limits_{j\in \Lambda}  \sum\limits_{i\in \bar{\Lambda}} \nabla h_j(x,\xi_i)^Tt^k  \rangle|$ for each $k,$ then it is suffices to prove that $\sum\limits^{\infty}_{k=0}\gamma_k<\infty.$ 
Now by Step~\ref{st4} and inequality (\ref{amj1}), it is established that
\begin{align*}
	H_j(x^{k+1})&\leq H_j(x^k)+\alpha \eta H^*_j(x^k,t^k),~~\forall~j\in \Lambda \\
	&\leq H_j(x^k)+\alpha_k \eta H^*_j(x^k,t^k)\\
	&\leq H_j(x^k)+\alpha_k \eta \sum\limits_{j\in \Lambda}  \sum\limits_{i\in \bar{\Lambda}} \nabla h_j(x^k,\xi_i)^Tt^k \\
	&\leq H_j(x^k)+\alpha_k \eta \langle \delta_{ij}, \; \sum\limits_{j\in \Lambda}  \sum\limits_{i\in \bar{\Lambda}} \nabla h_j(x,\xi_i)^Tt^k  \rangle, 
\end{align*}
therefore, \begin{equation}\label{3.30}
	H_j(x^{k})-H_j(x^{k+1})\geq-\alpha_k \eta \langle \delta_{ij}, \; \sum\limits_{j\in \Lambda}  \sum\limits_{i\in \bar{\Lambda}} \nabla h_j(x,\xi_i)^Tt^k  \rangle. 
\end{equation}
Since $x^k$ is a non Pareto critical point then $\nabla h_j(x,\xi_i)^Tt^k<0$ for each $i,j.$ Therefore, $$\alpha_k \langle \delta_{ij}, \; \sum\limits_{j\in \Lambda}  \sum\limits_{i\in \bar{\Lambda}} \nabla h_j(x,\xi_i)^Tt^k  \rangle<0.$$ Now, it can be taken as 
$$-\alpha_k \langle \delta_{ij}, \; \sum\limits_{j\in \Lambda}  \sum\limits_{i\in \bar{\Lambda}} \nabla h_j(x,\xi_i)^Tt^k  \rangle=\alpha_k |\langle \delta_{ij}, \; \sum\limits_{j\in \Lambda}  \sum\limits_{i\in \bar{\Lambda}} \nabla h_j(x,\xi_i)^Tt^k  \rangle|.$$ By using the latest expression and inequality (\ref{3.30}), it is obtained that 
\begin{equation}\label{3.31}
	H_j(x^{k})-H_j(x^{k+1})\geq\alpha_k \eta |\langle \delta_{ij}, \; \sum\limits_{j\in \Lambda}  \sum\limits_{i\in \bar{\Lambda}} \nabla h_j(x,\xi_i)^Tt^k  \rangle|. 
\end{equation}
By equation (\ref{3.31}) and $\gamma_k$ we get
$\gamma_k\leq\frac{2\beta}{\eta} \sum\limits^m_{j=1}(H_j(x^{k})-H_j(x^{k+1}) ).$  On adding from 1 to $M$ in the above inequality, where $M$ is any positive integer, it follows that
$$\sum\limits^M_{k=0}\gamma_k\leq\frac{2\beta}{\eta} \sum\limits^m_{j=1}(H_j(x^{0})-H_j(x^{M+1}) ).$$ 
\end{proof}
Since $\beta,\eta >0,$ $M$ is any positive integer, and $H(x')\leq H(x^k)$ for all k, then it is concluded that $\sum\limits^\infty_{k=0}\gamma_k<\infty.$ As $x'$ was an arbitrary element of $S$ therefore, $\{x^k\}$ quasi-Fejer converges to $S.$ From the discussion preceding this theorem, it is established that if ${x^k}$ is quasi-Fejer convergent to a set $S\subset\mathbb{R}^n$, then ${x^k}$ is bounded. Since $D$ is closed, and $\{x^k\}$ is bounded with $x^k\in D$ for all $k$ then there must exists an accumulation point, say $x^*\in D$ which will be the Pareto critical point of $H$ by Theorem \ref{th3.3}. Also by the convexity $x^*$ will be the weak Pareto optimal solution of $H$.
\section{Numerical results}\label{secnum}
In this section, numerical experiments are presented to evaluate the effectiveness of the proposed method and compare it to the weighted sum method. The focus is on finding the solution of $OWC_{P(U)}$ that corresponds to the solution of $P(U),$ using the projected gradient method. The goal is to verify how efficiently this method finds the Pareto front compared to the weighted sum method. To achieve this goal, the following methods were considered in the tests:\\
$	\checkmark$ To find the projected gradient descent direction at current iterative point $x^k$, it is necessary to solve the subproblem $\text{S}_P(x^k),$ which can be equivalently addressed using the following problems:
 \begin{equation}\label{eq32**}
	\underset{d\in \mathbb{R}^n}{\mathrm{min}}  \left(\beta \Theta_{x^k}(d)+\tfrac{1}{2}\|d\|^2 \right).
\end{equation}
Also, subproblem \eqref{eq32**} can be convert into the following equivalent form: $$P^{GD}(x^k):~~~~\displaystyle \min_{(r,~d)\in \mathbb{R}^{n+1}}\beta r+\frac{\|d\|^2}{2}$$
$$s.t.~~~~~~ \beta\big(h_j(x^k,\xi_i)+\nabla  h_j(x^k,\xi_i)^Td- H_{j}(x^k)\big)- \beta r\leq 0, ~\forall i\in \bar{\Lambda}, ~ \forall j\in \Lambda.$$ 	
Lagrangian function for $P^{GD}(x^k)$ is
\begin{align*}
	L(r,d,\lambda) &= \beta r+\frac{\|d\|^2}{2} + \sum\limits_{j\in \Lambda}  \sum\limits_{i\in \bar{\Lambda}} \lambda_{ij} \beta\bigg( h_j(x^k,\xi_i)+\nabla  h_j(x^k,\xi_i)^Td- H_{j}(x^k)- r\bigg).
\end{align*}
By using the KKT optimality conditions, it is obtained that
\begin{equation} \label{2}
	\sum\limits_{j\in \Lambda}  \sum\limits_{i\in\bar{\Lambda}} \lambda_{ij}=1,
\end{equation}
\begin{equation}\label{3}
\beta	\sum\limits_{j\in \Lambda}  \sum\limits_{i\in\bar{\Lambda}} \lambda_{ij}  \nabla  h_j(x^k,\xi_i)+ d=0,
\end{equation}
\begin{equation}\label{3*}
	\lambda_{ij}\geq0,~ \beta\big(h_j(x^k,\xi_i)+\nabla  h_j(x^k,\xi_i)^Td- H_{j}(x^k)- r\big)\leq 0,~  ~\forall i\in \bar{\Lambda},~ \forall j\in \Lambda,
\end{equation}
\begin{equation}\label{3**}
	\lambda_{ij}\beta\bigg( h_j(x^k,\xi_i)+\nabla  h_j(x^k,\xi_i)^Td- H_{j}(x^k)- r\bigg)=0, ~\forall i\in \bar{\Gamma}, ~ \forall j\in \Gamma.
\end{equation}
Problem $P^{GD}(x^k)$ has a unique solution, $(r(x^k),d(x^k)).$ As this is a convex problem and has a Slater point, there exists a KKT multiplier $\lambda=\lambda(x^k)=\lambda_{ij}(x^k),$ together $r=r(x^k)$ and $d=d(x^k),$ satisfies conditions (\ref{2}), (\ref{3}), (\ref{3*}), and (\ref{3**}). In particular, by using (\ref{3}), it can be obtained that
\begin{equation}\label{newt}
	d(x^k)=- \beta\sum\limits_{j\in \Lambda}  \sum\limits_{i\in\bar{\Lambda}} \lambda^k_{ij}(x^k)  \nabla  h_j(x^k,\xi_i).
	\end{equation}
	Two cases may arise:
	
	\textbf{Case 1:} If $d(x^k)$ lies in $D$, then $t(x^k) = d(x^k)$.
	
\textbf{Case 2:} If $d(x^k)$ does not lie in $D$, then to find $t(x^k)$, the projection of $d(x^k)$ onto the set $D$ is taken as $P_D(d(x^k)).$
	Equivalently, $t(x^k)$ is obtained as:
	\[t(x^k) = \text{Proj}_{D} \left( - \beta \sum\limits_{j \in \Lambda} \sum\limits_{i \in \bar{\Lambda}} \lambda_{ij} \nabla h_j(x^k, \xi_i) \right).\]
In the upcoming computation, for any convex set $D,$ $\text{Proj}_D(d(x^k)) = \underset{y \in D}{\mathrm{argmin}} \|y - d(x^k)\|^2_2$ is used as the projection of $d^k$ on $D.$\\
$	\checkmark$ If stopping criterion is not satisfied at $x^k$ $(\|t^k\|\not<10^{-4})$ i.e., $x^k$ is not a approximate critical point, then to find the step length size at $x^k$ in projected gradient descent direction, the following condition is used:
\begin{equation}\label{13at k}
	H_j(x^k+\alpha t^k)\leq  H_j(x^k)+  \alpha \eta H^*_j(x^k,t).
\end{equation}
According to this equation $\alpha_k$ will be the largest $\alpha \in \{ \frac{1}{2^r} : r=1,2,3,\ldots\} $ satisfying inequality (\ref{13at k}).\\
	$\checkmark$ To address the MOP using the weighted sum method, the following single objective optimization problem is considered:
$$\min_{x \in D} \left( w_1 H_1(x) + w_2 H_2(x) + \dots + w_m H_m(x) \right),$$
where the weights $(w_1, w_2, \dots, w_m)$ are non-negative and at least one weight is non-zero. The technique developed by \cite{bazaraa1982algorithm} is used with the initial guess $x^0 = \frac{1}{2}(lb + ub)$. For bi-objective optimization problems, weights $(1,0)$, $(0,1)$, and 98 random weights uniformly distributed within the square $[0,1] \times [0,1]$ are used. Similarly, for three-objective optimization problems, four specific weights $(1,0,0)$, $(0,1,0)$, $(0,0,1)$, and 97 random weights uniformly distributed within the cube $[0,1] \times [0,1] \times [0,1]$ are considered.\\
	$\checkmark$Using the above information, a comparison is presented between the projected gradient method and the weighted sum method through numerical examples. Python code is written to solve these examples using Algorithm \ref{algo1}. The stopping criterion for the calculations is either a tolerance of $\epsilon=10^{-4}$ or a maximum of 5000 iterations. Therefore, the algorithm stops when either $\|t(x^k)\| < 10^{-4}$ or the maximum of 5000 iterations is reached. It is understood that the solution to the MOP is not a single point but a set of Pareto optimal solutions. To obtain a well-distributed approximation of the Pareto front, a multi-start technique is used. Specifically, 100 uniformly distributed random points between $lb$ and $ub$ (where $lb, ub \in D\subset \mathbb{R}^n$ and $lb < ub$) are generated, and Algorithm \ref{algo1} is executed separately at each point. The approximate Pareto front is the set of nondominated points that are not dominated. Now by using the following examples, the approximate Pareto fronts produced by Algorithm \ref{algo1} are compared with those produced by the weighted sum method.

\begin{example}\label{ex1}
	Consider $H(x)=(H_1(x),H_2(x)),$ where
	$H_1(x)= \displaystyle\max_{\xi_i\in U}h_1(x,\xi_i),$ $H_2(x)=\displaystyle\max_{\xi_i\in U}h_2(x,\xi_i)$ and $h_1(x,\xi_i)=(x_1-\xi^1_i)^2+(x_2-\xi^2_i)^2,$ $h_2(x,\xi_i)=\xi^1_ix_1^2+\xi^2_ix^2,$ and $U=\{\xi_1=(1,2),~\xi_2=(1,1)\}$ with $\xi_i=(\xi^1_i,\xi^2_i),$ $i=1,2.$
	$$OWC_{P(U)}:~~\min ~H(x)$$
	$$~~~~~\text{s. t.}~~~~~~ ~-5\leq x_1,~x_2\leq 10.$$
%
		Consider $x^0=(-4.4,4.4)^T$. Then, $h_{1}(x^0,\xi_1)=23.5$, $h_{1}(x^0,\xi_2)=58.5$, $h_{2}(x^0,\xi^{1})=59$, and $h_{2}(x^0,\xi_2)=59$. \\
	Hence, $H(x^0)=(58.5,59)^T$. One can observe that \[I_1(x^0)=\{2\} \text{ and } I_2(x^0)=\{1,2\}.\]
%
%
%
%
%
%
 	Solution of $P^{GD}(x^0)$ is obtained as $d^0=(11.871, -2.287)^T$, $r^0=-140.123$. 
		According to the problem $D = \{(x_1, x_2) : -5 \leq x_1, x_2 \leq 10\},$ and $d^0=(11.871, -2.287)\not\in D=\{(x_1,x_2)\in\mathbb{R}^2:-5\leq x_1,~x_2\leq 10\}.$ Therefore, it is required to find the value of $\text{Proj}_D(11.871, -2.287),$ which is the projection of the point $d^0$ onto the set $D.$ The projection of the point $d^0$ onto the set $D$ can be determined as follows:
	
	\[
	\text{Proj}_D(d^0) = \arg \min_{y \in D} \| y - d^0 \|_2^2= (10, -2.287).
	\]
%
%
%
%
%
	 Consequently, $t^0=(10, -2.287)$. Also, it can be find that $\alpha_0=0.5$ satisfies (\ref{13}). Hence, next iterating point is $x^1=x^0+\alpha_0 t^0=(0.6,3.2565)^T$. \\\\
	One can observe that $H(x^1)=( 9.7023,21.5447)^T<H(x^0)$. Using the stopping criterion $\|t^k\|<10^{-4}$, the final solution is obtained as $x^*=(1.9673,1.9673)^T$ after $5$ iterations.\\
	
	In the next, it is demonstrated that $x^*=(1.9673,1.9673)^T$ is a week efficient solution of this problem. Observe that $h_{1}(x^*,\xi_{1})=3.0353=h_{1}(x^*,\xi_{2})$ and $h_{2}(x^*,\xi_{1})=15.4957=h_{2}(x^*,\xi_{2})$. These imply $I_1(x^*)=\{1,2\}$ and $I_2(x^*)=\{1,2\}$. Hence,
	\[\partial H_1(x^*)=conv \{(2.6343,-3.4653)^T, (-3.4653,2.6343)^T\}\] and \[\partial H_2(x^*)=conv \{(4.1345,13.4092)^T, (13.4092,4.1345)^T\}.\] Clearly,
	$$e_1=0.5(2.6343,-3.4653)^T+0.5(-3.4653,2.6343)^T=(-0.4155,-0.4155)^T\in\partial H_1(x^*).$$
	 Similarly, $e_2=(8.5691,8.5691)^T\in \partial H_2(x^*)$. Choose $\lambda_1=\frac{85691}{89846}$. Then, $0<\lambda_1<1$, $\lambda_1e_1+(1-\lambda_1)e_2=0$. This implies $$0\in conv \{\partial H_1(x^*),\partial H_2(x^*)\}$$ and hence $x^*$ is a weak efficient solution.
	 \begin{figure}[!htbp]
	 	\centering
	 		\includegraphics[width=3.0in,height=1.4in]{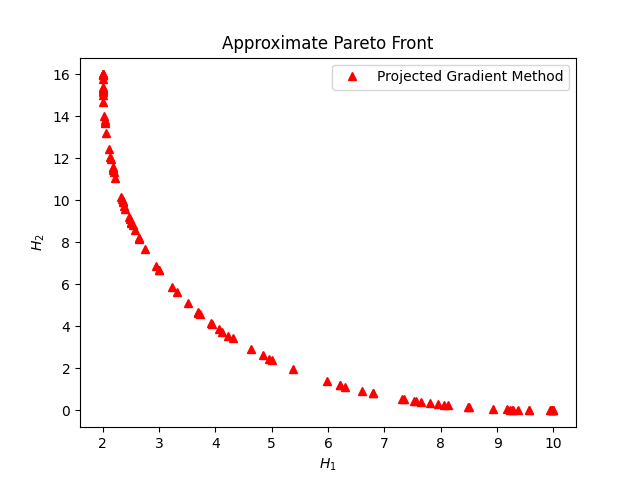}
	 	\includegraphics[width=3.0in,height=1.4in]{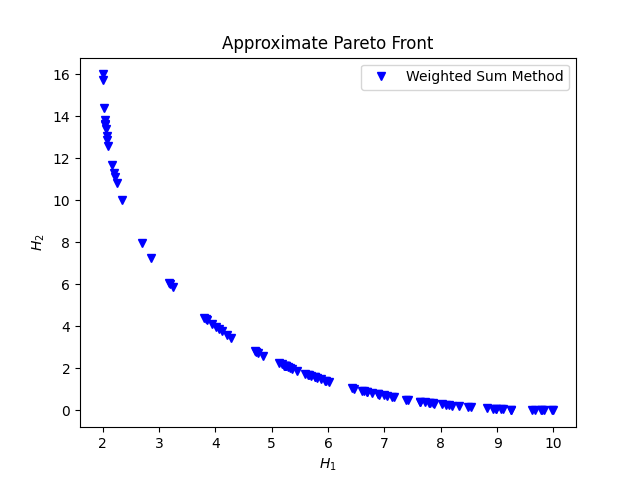}
	 	\caption{\textbf{Comparison of approximate Pareto fronts generated by Algorithm \ref{algo1} and weighted sum method for Example \ref{ex1}.}}
	 	\label{fig1}
	 \end{figure}
	
	\end{example}
\begin{example}\label{ex2} Consider
		$H(x)=(H_1(x),H_2(x)),$ where
		$H_1(x)= \displaystyle\max_{\xi_i\in U}h_1(x,\xi_i),$ $H_2(x)=\displaystyle\max_{\xi_i\in U}h_2(x,\xi_i),$   and $h_1(x,\xi_i)=(x-\xi_i)^2,$ $h_2(x,\xi_i)=-x^2-\xi_i x,$ and $U=\{\xi_1=-5,~\xi_2=2\}.$
			$$OWC_{P(U)}:~~\min ~H(x)$$
		$$~~~~~\text{s. t.}~~~~~~ ~-9\leq x\leq 5.$$
		 \begin{figure}[!htbp]
			\centering
			\includegraphics[width=3.0in,height=1.4in]{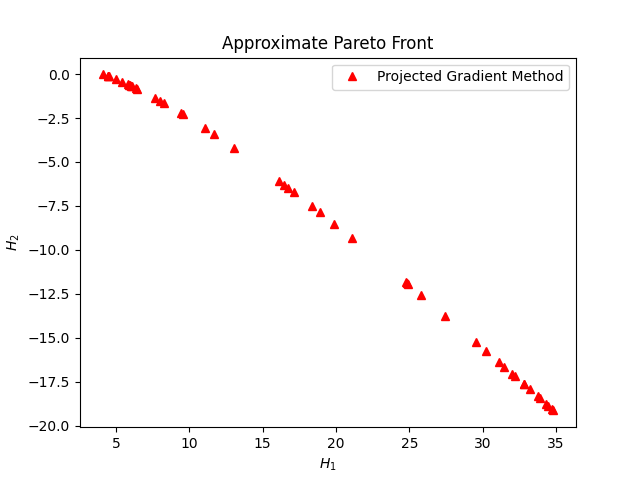}
			\includegraphics[width=3.0in,height=1.4in]{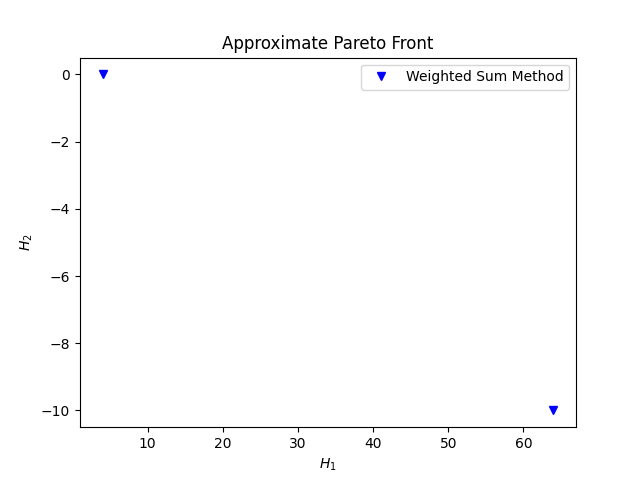}
			\caption{\textbf{Comparison of approximate Pareto fronts generated by Algorithm \ref{algo1} and weighted sum method for Example \ref{ex2}.}}
			\label{fig2}
		\end{figure}
\end{example}
\begin{example}\label{ex3}
	Consider $H(x)=(H_1(x),H_2(x),H_3(x)),$ where
	$H_1(x)= \displaystyle\max_{\xi_i\in U}h_1(x,\xi_i),$ $H_2(x)=\displaystyle\max_{\xi_i\in U}h_2(x,\xi_i),$  $H_3(x)=\displaystyle\max_{\xi_i\in U}h_3(x,\xi_i),$ and $h_1(x,\xi_i)= x_1^2+\xi^1_ix_2^4+\xi^1_i\xi^2_ix_1x_2,$ $h_2(x,\xi_i)=5x_1^2+\xi^1_ix_2^2+\xi^2_ix_1^4x_2,$ 	$h_3(x,\xi_i)=e^{-\xi^1_ix_1+\xi^2_ix_2}+x_1^2-\xi^1_ix_2^2$ and $\xi_i\in U=\{(2,3),~(4,5),~(2,0)\}$ with $\xi_i=(\xi^1_i,\xi^2_i),$ $i=1,~2,~3.$
		$$OWC_{P(U)}:~~\min ~H(x)$$
	$$~~~~~\text{s. t.}~~~~~~ ~-11\leq x_1,x_2\leq 5.$$
%
	 \begin{figure}[!htbp]
		\centering
		\includegraphics[width=3.2in,height=1.47in]{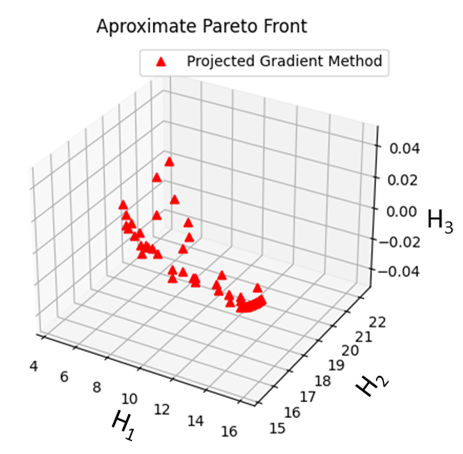}
		\includegraphics[width=3.2in,height=1.47in]{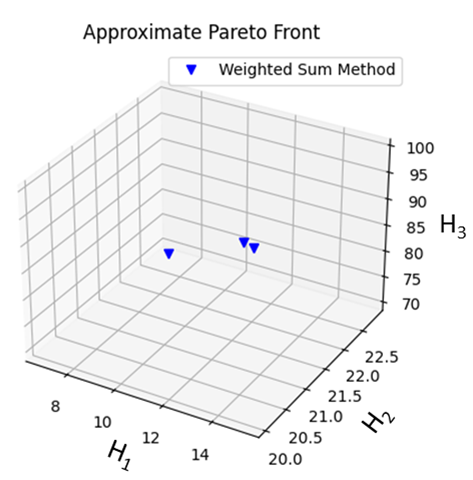}
		\caption{\textbf{Comparison of approximate Pareto fronts generated by Algorithm \ref{algo1} and weighted sum method for Example \ref{ex3}.}}
		\label{fig3}
	\end{figure}
\end{example}
\begin{example}\label{ex4}
	Consider $H(x)=(H_1(x),H_2(x),H_3(x)),$ where
	$H_1(x)= \displaystyle\max_{\xi_i\in U}h_1(x,\xi_i),$ $H_2(x)=\displaystyle\max_{\xi_i\in U}h_2(x,\xi_i),$  $H_3(x)=\displaystyle\max_{\xi_i\in U}h_3(x,\xi_i),$ and $h_1(x,\xi_i) =100\xi^1_i (x_2-x_1^2)^2+\xi^2_i(1-x_1)^2,$ $h_2(x,\xi_i)=(x_2-\xi^1_i)^2+\xi^2_ix_1^2,$ 	$h_3(x,\xi_i)=\xi^1_ix_1^2+3\xi^2_ix_2^2$ and $\xi_i\in U=\{(2,3),~(1,2),~(4,5)\}$ with $\xi_i=(\xi^1_i,\xi^2_i),$ $i=1,~2,~3.$
	$$OWC_{P(U)}:~~\min ~H(x)$$
	$$~~~~~\text{s. t.}~~~~~~ ~-10\leq x_1,x_2,x_3\leq 10.$$
	
	 \begin{figure}[!htbp]
		\centering
		\includegraphics[width=3.2in,height=1.47in]{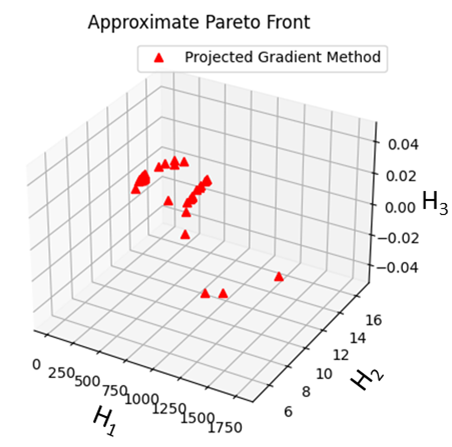}
		\includegraphics[width=3.2in,height=1.47in]{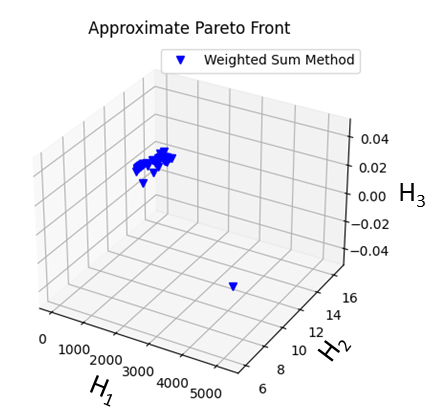}
		\caption{\textbf{Comparison of approximate Pareto fronts generated by Algorithm \ref{algo1} and weighted sum method for Example \ref{ex4}.}}
		\label{fig4}
	\end{figure}
	\end{example}
	From the numerical examples provided above, it is evident that Example \ref{ex1} and Example \ref{ex4} are convex, on the other hand Example \ref{ex2} and Example \ref{ex3} are nonconvex. By the figures \big(Figure~\ref{fig1}, Figure~\ref{fig2}, Figure~\ref{fig3}, and Figure~\ref{fig4}\big) generated by projected gradient method and weighted sum method corresponding to each example, it is clear that PGM successfully generates good approximate Pareto fronts in both convex and nonconvex case, and weighted sum method fails to generate good approximate Pareto fronts in non convex case \big(Figure~\ref{fig2}, Figure~\ref{fig3}\big). Algorithm \ref{algo1} offers a distinct advantage by eliminating the challenging task of selecting pre-order weights, a requirement in the weighted sum method. Consequently, without the burden of pre-order weight selection, Algorithm \ref{algo1} consistently produces an efficient approximations of Pareto fronts. The primary objective of identifying the approximate Pareto front for the uncertain multiobjective optimization problem \(P(U)\) with a finite parameter uncertainty set has been achieved through the utilization of \(OWC_{P(U)}\) and PGM.
	\section{Conclusions}\label{s4}
In this article, authors tackled a constrained uncertain multiobjective optimization problem, \( P(U) \), with convex constraints using a robust optimization approach and the projected gradient method (PGM). The authors demonstrated that solving the objective-wise worst-case cost type robust counterpart, \( OWC_{P(U)} \), suffices to solve \( P(U) \). By focusing on \( OWC_{P(U)} \) with PGM, authors presented a novel subproblem for the projected descent direction and an Armijo-type inexact line search for step length selection. With the help of projected descent direction and step length size a projected descent algorithm is written. It is shown that the sequence generated by the projected gradient descent algorithm converges to a weak Pareto optimal solution of $OWC_{P(U)}$ as well as the robust weak Pareto optimal solution of $P(U).$ The global convergence is also proved under convexity assumptions. Proposed method is also verified some numerical examples, and made a comaparison with existing weighted sum method. By numerical tests, it is also observed that projected gradient method give an efficient Pareto fronts in both convex and nonconvex problem. On the other hand, weighted sum method for $P(U)$ fails to generate good approximate for nonconvex problem. Consequently, PGM for ${P(U)}$ does not require predetermined weighting factors or any other kind of predetermined ranking or ordering information for objective functions, eliminating another drawback of scalarization methods. 
\par In this study, the Armijo-type inexact line search technique is used to determine the step length size selection. Future research will focus on developing Wolfe and Zoutendijk conditions for step length size determination.

Also, in the presented study, the focus was on solving an uncertain multiobjective optimization problem under a finite uncertainty set. Addressing the solution of such problems under an infinite uncertainty set remains a subject for future investigation.
\section*{Availability of data and materials}
Not applicable.
\section*{Acknowledgment:}
This research work of the first author is supported by Govt. of India CSIR fellowship, Program No. 09/1174(0006)/2019-EMR-I, New Delhi, India.
\section*{Conflict of interest:}
The authors declare no conflicts of interest.
\bibliographystyle{elsarticle-num}
\bibliography{projectedgradient_reference}

\end{document}